\documentclass[12pt,twoside,reqno]{amsart}
\usepackage[colorlinks=true,citecolor=blue]{hyperref}
\usepackage{mathptmx, amsmath, amssymb, amsfonts, amsthm, mathptmx, enumerate, color}
\usepackage{graphicx}
\usepackage{epstopdf}

\newtheorem{theorem}{Theorem}[section]

\newtheorem{lemma}{Lemma}[section]
\newtheorem{proposition}{Proposition}[section]
\theoremstyle{definition}
\newtheorem{definition}{Definition}[section]
\newtheorem{example}{Example}[section]
\newtheorem{remark}{Remark}[section]

\numberwithin{equation}{section}

\usepackage{bbm}
\usepackage{tikz}
\usepackage{pgfplots,pgfplotstable}
\usetikzlibrary{arrows,positioning,chains,fit,shapes,calc,decorations}
\usepackage{subfig}

\usepackage{algorithm,algpseudocode}
\algnewcommand\algorithmicinput{\textbf{INPUT:}}
\algnewcommand\INPUT{\item[\algorithmicinput]}
\algrenewcommand\Return{\State \algorithmicreturn{} }%

\algtext*{EndWhile}
\algtext*{EndIf}
\algtext*{EndFor}

\newcommand{\RR}{\mathbb{R}}
\newcommand{\NN}{\mathbb{N}}  
\newcommand{\NNN}{\overline{\NN}}  
\def\D{{\mathcal I}}
\def\Pcal{{\mathcal P}}
\def\Fcal{{\mathcal F}}

\pgfplotsset{compat=1.16}

\begin{document}
\setcounter{page}{1}

\vspace*{1.0cm}
\title[Optimal transports bounds for fixed point iterations]{Universal bounds for fixed point iterations via  optimal transport metrics}

\author[M. Bravo, T. Champion, R. Cominetti]{ Mario Bravo$^{1}$, Thierry Champion$^2$, Roberto Cominetti$^{1,*}$} 
\maketitle
\vspace*{-0.6cm}

\begin{center}
{\footnotesize {\it
$^1$Facultad de Ingenier\'ia y Ciencias, Universidad Adolfo Ib\'a\~nez, Diagonal Las Torres 2640, Santiago, Chile.
 \\
$^2$Laboratoire Imath, U.F.R. des Sciences et Techniques, Universit\'e de Toulon, Av. de l'Universit\'e, B.P. 20132, 83957 La Garde cedex, France.
}}\end{center}

\vskip 3mm
\centerline{\em  Dedicated to the memory of Professor Ronald E. Bruck}

\vskip 4mm 
{\small\noindent {\bf Abstract.}
We present a self-contained analysis of a particular family 
of metrics over the set of non-negative integers. We show that these metrics, which are 
defined through a nested sequence of  optimal transport problems, provide tight estimates for
 general Krasnosel'skii-Mann fixed point iterations for non-expansive maps. We also describe 
 some of their very special properties, including their monotonicity and the  so-called {\em convex 
 quadrangle inequality} that yields a greedy algorithm for computing them efficiently. 
 
 \vskip 2mm 
\noindent {\bf Keywords.} Fixed-point iterations, non-expansive maps, error bounds, convergence rates, optimal transport metrics.

}

\renewcommand{\thefootnote}{}
\footnotetext{ 
$^*$Corresponding author.
\par
E-mail addresses: 
{\scriptsize 
\tt{<mario.bravo@uai.cl>}, 
\tt{<champion@univ-tln.fr>}, 
\tt{<roberto.cominetti@uai.cl>} 
}
\par
Received xxxxx  yy, 2021; Accepted  xxxxx yy, 2021. 
}

\section{Introduction}

This paper studies a special family of metrics over the set of non-negative integers, in which the distances $d_{m,n}$ between 
$m,n\in\NN$ are defined recursively through a nested family of optimal transport problems. 
These metrics were first introduced in a remarkable paper by Baillon \& Bruck \cite{bb92} with the aim
of establishing convergence rates for the Krasnosel'skii fixed-point iteration, 
although their metric nature and the connection with optimal transport were not noticed at that time. 
The metric properties of the $d_{m,n}$'s  were studied extensively in Aygen-Satik's thesis \cite{aig}, with a {\em tour-de-force} that required very 
long and highly technical proofs. Several of these results were revisited in \cite[Bravo \& Cominetti]{brc} with a much simpler approach based on
optimal transport.
The goal of this paper is to complete this program by presenting a full and self-contained analysis of these metrics,
and their implications for fixed point iterations. 
In doing so, we expand the scope of \cite{bb92,aig,brc} by considering a much larger family of 
iterations and metrics.

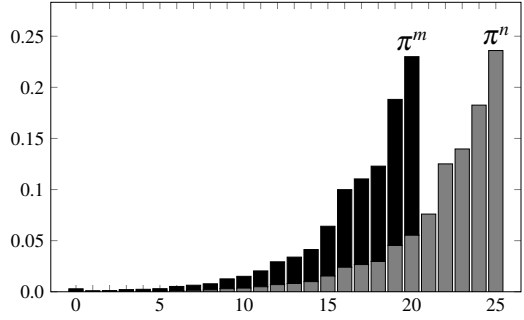
\begin{figure}[ht]
\centering
\begin{tikzpicture}[scale=0.6]
    \begin{axis}[ybar stacked,
        ymin=0,
        xmin=-1.5,
        xmax= 26.5,
        width  = 12cm,
        height = 8cm,
        bar width=9pt,
        xtick = data,
        table/header=false,
        table/row sep=\\,
        xticklabels from table={ 0\\ \  \\ \  \\ \ \\ \ \\ 5\\ \ \\ \ \\ \ \\ \ \\ 10\\ \ \\ \ \\ \ \\ \ \\ 15\\ \ \\ \ \\ \ \\ \ \\20\\ \ \\ \ \\ \ \\ \ \\25\\}{[index]0},
        yticklabels from table={ 0.0\\ 0.0 \\ 0.05\\ 0.1\\ 0.15\\0.2 \\ 0.25\\ 0.3\\}{[index]0},
        enlarge y limits={value=0.2,upper},
    ]
 \node[anchor=south] at (axis cs:20,0.23) {\Large $\pi^m$};
 \node[anchor=south] at (axis cs:25,0.235) {\Large $\pi^n$};

    \addplot[fill=gray] table[x expr=\coordindex,y index=0]{ 
    0.0007\\
    0.0002\\
    0.0003\\
    0.0005\\
    0.0006\\
    0.0007\\
    0.0012\\
    0.0015\\
    0.0019\\
    0.0030\\
    0.0036\\
    0.0049\\
    0.0070\\
    0.0081\\
    0.0099\\
    0.0154\\
    0.0240\\
    0.0266\\
    0.0296\\
    0.0453\\
    0.0553\\
    0.0760\\
    0.1251\\
    0.1397\\
    0.1826\\
    0.2360\\
};
    \addplot[fill=black] table[x expr=\coordindex,y index=0]{
  0.0022\\
   0.0008\\
   0.0008\\
   0.0016\\
   0.0018\\
   0.0023\\
   0.0040\\
   0.0049\\
   0.0059\\
   0.0096\\
   0.0115\\
   0.0155\\
   0.0223\\
   0.0258\\
   0.0314\\
   0.0487\\
   0.0759\\
   0.0838\\
   0.0933\\
   0.1429\\
   0.1747\\
    0\\
    0\\
    0\\
    0\\
    0\\
};

    \end{axis}
\end{tikzpicture}
\vspace{-2ex}
\caption{\label{pis}The distributions $\pi^m$ and $\pi^n$ ($m=20, n=25$).}
\end{figure}

To get into the matter, throughout this paper we consider a fixed sequence $(\pi^n)_{n\in\NN}$ where 
each $\pi^n=(\pi_i^n)_{i\in\NN}$ is a discrete probability distribution on the set of non-negative integers,
with support included in $\{0,\ldots,n\}$ and with $\pi^n\neq\pi^m$ for $m\neq n$. Some results  
assume in addition that these distributions  gradually drift  their mass  towards larger integers,
namely (see Figure \ref{pis})
$$(\forall n\geq 1)\quad \mbox{$\pi_n^n>0$ and $0\leq \pi_i^{n}\leq\pi_i^{n-1}$ for $i<n$.}\leqno\mbox{\sc(h)}$$

\subsection{Recursive optimal transports}
Starting with $d_{-1,j}=d_{j,-1}= 1$ for  $j\in\NN$ and $d_{-1,-1}=0$, we consider the double-indexed family of reals 
$d_{m,n}$ defined recursively for $m,n\in\NN$ through the optimal transport problems 
\begin{equation*}\label{dmn}
d_{m,n}=\min_{z\in\Fcal_{m,n}} C_{m,n}(z)\triangleq \sum_{i=0}^m\sum_{j=0}^n z_{i,j}\,d_{i-1,j-1}\leqno (\Pcal_{m,n})
\end{equation*}
where $\Fcal_{m,n}$ is the polytope of transport plans sending
$\pi^m$ to $\pi^n$, that is, the set of all $z=(z_{i,j})_{i=0,\ldots,m; j=0,\ldots n}$
such that $z_{i,j}\geq 0$ and
$$\begin{array}{ll}
\sum_{j=0}^n z_{i,j}=\pi_i^m&\mbox{for all } i=0,\ldots,m;\\[1.5ex]
\sum_{i=0}^m z_{i,j}=\pi_j^n&\mbox{for all } j=0,\ldots,n.
\end{array}
$$
\begin{figure}[h]
\centering
\begin{tikzpicture}[thick, -,shorten >= 1pt,shorten <= 1pt,scale=0.65,every node/.style={scale=0.7}]
 \node[draw,circle] (i0) at (0,0) [label=left: {\Large $\pi_0^m$}] {$0$} ;
 \node[] (i1) at (0,-0.9) {$\vdots$};
 \node[draw,circle] (im) at (0,-2) [label=left: {\Large $\pi_m^m$}] {$m$};

 \node[draw,circle] (j0) at (4,0) [label=right: {\Large $\pi_0^n$}] {$0$};
 \node[] (j1) at (4,-0.9) {$\vdots$};
 \node[draw,circle] (jm) at (4,-2) [label=right: {\Large $\pi_j^n$}] {$j$};
 \node[] (j2) at (4,-2.9) {$\vdots$};
 \node[draw,circle] (jn) at (4,-4) [label=right: {\Large $\pi_n^n$}] {$n$};

\draw[->] (i0) -- (j0);
\draw[->] (i0) -- (jm);
\draw[->] (i0) -- (jn);
\draw[->] (im) -- (j0);
\draw[->] (im) -- (jm);
\draw[->] (im) -- (jn);
\end{tikzpicture}
\caption{\label{feasible_plan}The optimal transports for $d_{m,n}$.}
\end{figure}
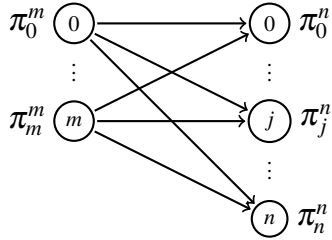

Each sequence $(\pi^n)_{n\in\NN}$ produces a different set of $d_{m,n}$'s. 
In particular we have $d_{0,n}=1 - \pi_0^n$, though in general
the other $d_{m,n}$'s do not admit a simple explicit formula. 
We observe that a simple induction yields the symmetry $d_{m,n}=d_{n,m}$
so it suffices to compute $d_{m,n}$ for $m\leq n$. Also, since the transports
$z_{i,j}$ add up to one, the cost $C_{m,n}(z)$ is a convex combination of the 
previous $d_{i-1,j-1}$'s and inductively we get $d_{m,n}\in[0,1]$ with $d_{n,n}=0$.

\subsection{Fixed point iterations.}\label{secKM} 
The optimal transports $d_{m,n}$ arise in connection with fixed point iterations. 
Namely, let  $T:C\to C$ be a non-expansive map on a bounded convex domain $C\subseteq X$ of a normed space $(X,\|\cdot\|)$. 
Notice that when $T:X\to X$ is defined on the full space and has some fixed point $x^*=Tx^*$, one may take $C=B(x^*,r)$ as any ball  
centered at $x^*$ with radius $r\geq 0$. On the other hand, by rescaling the norm by a factor $1/\mathop{\rm diam}(C)$ we may assume without 
loss of generality that $\mathop{\rm diam}(C)=1$, which we do from now on.

Starting from $x^0,y^0\in C$, and setting by convention $Tx^{-1}=y^0$, the general Krasnosel'skii-Mann fixed point iteration (see \cite[Mann]{man} and \cite[Krasnosel'skii]{kra})
recursively builds a sequence $x^n$  as a convex combination of the images of the previous iterates
$$\mbox{$x^n=\sum_{i=0}^n\pi_i^n\,Tx^{i-1}$}\qquad\forall n\geq 1.\leqno\mbox{\sc (km)}$$

The following straightforward result brings forward the connection with the recursive optimal 
transports, showing how they provide bounds 
for the distance between the {\sc (km)} iterates and for the fixed point residuals.
\begin{theorem} \label{prop0}
For all $m,n\in\NN$ we have $\|x^m-x^n\|\leq d_{m,n}$  and
\begin{equation}\label{residual}
\|x^n-Tx^n\|\leq R_n\triangleq \mbox{$\sum_{i=0}^n\pi_i^nd_{i-1,n}$}.
\end{equation}
\end{theorem}
\begin{proof} 
Let us first show that  $\|x^m-x^n\|\leq d_{m,n}$. The cases $m=0$ and $n=0$ follow 
directly from $\|x^{0}-x^{n}\|\leq \mathop{\rm diam}(C)=1=d_{-1,n}$. For the general case, we 
observe that each transport plan $z$ from $\pi^m$ to $\pi^n$ yields the estimate
\begin{eqnarray*}
\|x^m-x^n\|&=&\mbox{$\|\sum_{i=0}^m\pi^m_iTx^{i-1}- \sum_{j=0}^n\pi^n_jTx^{j-1}\|$}\\
&=&\mbox{$\|\sum_{i=0}^m\sum_{j=0}^nz_{i,j}(Tx^{i-1}-Tx^{j-1})\|$}\\
&\leq&\mbox{$\sum_{i=0}^m\sum_{j=0}^nz_{i,j}\|Tx^{i-1}-Tx^{j-1}\|$}\\
&\leq&\mbox{$\sum_{i=0}^m\sum_{j=0}^nz_{i,j}d_{i-1,j-1}$}
\end{eqnarray*}
where the last inequality uses the non-expansivity of $T$  and assumes inductively that we already have $\|x^{i-1}-x^{j-1}\|\leq d_{i-1,j-1}$
for the previous iterates (for $i=0$ or $j=0$ use the convention $Tx^{-1}=y^0$ and the coarse
estimate $\|y^{0}-Tx^{k}\|\leq \mathop{\rm diam}(C)=1=d_{-1,k-1}$).
Minimizing over $z$ we get $\|x^m-x^n\|\leq d_{m,n}$ and then the proof is completed
by induction. The bound \eqref{residual}  follows directly from the triangle inequality and non-expansivity, namely
$$\|x^n-Tx^n\|=\mbox{$\|\sum_{i=0}^n\pi^n_i(Tx^{i-1}-Tx^n)\|$}\leq\mbox{$\sum_{i=0}^n\pi_i^nd_{i-1,n}$}.$$ 

\vspace{-3ex}
\end{proof}

We emphasize that the bounds $d_{m,n}$ and $R_n$ are universal in the sense that they only depend on the sequence 
$(\pi^n)_{n\in\NN}$ and not on the particular map $T$ being considered. Moreover, the estimate $\|x^m-x^n\|\leq d_{m,n}$
turns out to be tight and cannot be improved unless we restrict the class of maps or spaces. 
In addition, we will show that under {\sc (h)} the bound $\|x^n-Tx^n\|\leq R_n$
is also sharp, in which case $R_n$ captures the exact convergence rate of {\sc (km)}.

The recursion {\sc (km)} is very general and includes among others the Krasnosel'skii iteration (see \cite[Krasnosel'skii]{kra} and \cite[Bruck]{bru})
 $$x^{n}=(1\!-\!\alpha_n)x^{n-1}\!+\alpha_n Tx^{n-1}, \qquad\alpha_n\in(0,1),$$ 
as well as Halpern's method (see  \cite[Halpern]{hal} and \cite[L\'opez {\em et al.}]{lmh})
$$x^n=(1\!-\!\beta_n) y^0+\beta_nTx^{n-1},\qquad\beta_n\in(0,1),$$  
and the 2-step iteration of Ishikawa  (see \cite[Ishikawa]{ish0})
$$\left\{\begin{array}{rcl}
x^{2n+1}&=&\!(1\!-\!\beta_n)x^{2n}\!+\beta_n Tx^{2n}\\
x^{2n+2}&=&\!(1\!-\!\alpha_n)x^{2n}\!+\alpha_n Tx^{2n+1}
\end{array}\right.,\quad 0\leq\alpha_n\leq\beta_n\leq 1.
$$

 
The  Krasnosel'skii iteration is the case where $\pi^{n}=(1\!-\!\alpha_{n})\pi^{n-1}\!+\alpha_{n}\delta^{n}$ with $\delta^{n}$ a Dirac at $n$, 
for which we have explicitly $\pi_i^n= \alpha_i\prod_{k=i+1}^n (1-\alpha_k) $ with $\alpha_0=1$ and $\pi_n^n=\alpha_n$. 
Its convergence rate was studied in a series of papers \cite{bb92,aig,brc,bb96,csv}. 
A first result in \cite[Baillon \& Bruck]{bb92} for $\alpha_n\equiv \alpha$ showed that $\|x^n-Tx^n\|\sim O(\frac{1}{\ln n})$. 
Shortly after, in \cite[Baillon \& Bruck]{bb96} this rate was
improved to $\|x^n-Tx^n\|\sim O(\frac{1}{\sqrt{n}})$,  with an 
explicit upper bound recently extended  to non-constant $\alpha_n$'s  in \cite[Cominetti {\em et al.}]{csv}.
The papers \cite{bb92,bb96} are two gems among the many lasting contributions of Professor Ronald E. Bruck 
to the field of fixed point theory. A more 
detailed discussion of these previous results, as well as for Halpern iteration, is postponed to section \S\ref{s51}. 

\subsection{Summary of results}
This paper investigates the metric properties of the optimal transport bounds $d_{m,n}$.
Section \S\ref{smetric} shows that $(m,n)\mapsto d_{m,n}$ defines
a distance on $\NNN\triangleq \NN\cup\{-1\}$, which in turn leads to 
an alternative characterization via Kantorovich-Rubinstein duality.

This dual characterization is exploited in section \S\ref{sfin} to prove that the estimates  
$\|x^m-x^n\|\leq d_{m,n}$ in Theorem \ref{prop0} are the best possible, by 
building a particular non-expansive map $T$ and a corresponding {\sc (km)} 
sequence that attains these bounds with equality. Moreover, we prove that under {\sc (h)} 
the bounds $\|x^n-Tx^n\|\leq R_n$ for the residuals are also tight.

Section \S\ref{s51} illustrates these bounds in some specific iterations. 
For Halpern with $\beta_n=\frac{n}{n+1}$
we obtain $R_n=\frac{H_{n+1}}{n+1}\sim O(\frac{\ln n}{n})$ with $H_n=\sum_{k=1}^n\frac{1}{k}$ 
 the
$n$-th harmonic number, while for $\beta_n\!=\!\frac{n}{n+2}$ we have 
$R_n\!=\!\frac{4}{n+1}(1\!-\!\frac{H_{n+2}}{n+2})\!\sim\! O(\frac{1}{n})$,
both rates being tight for the corresponding $\beta_n$'s.

In section \S\ref{squadrangle} we show that under {\sc (h)}
the  $d_{m,n}$'s enjoy two relevant additional properties: the monotonicity as $m$ and $n$ get farther apart,  
and the so-called {\em convex quadrangle inequality} 
\begin{equation*}\label{4point}
d_{i,l}+d_{j,k}\leq d_{i,k}+d_{j,l}\qquad   \mbox{ for all $i\leq j\leq k\leq l$.}\leqno (Q)
\end{equation*}
Note that for $j=k$ this is just a triangle inequality.

\vspace{-1ex}
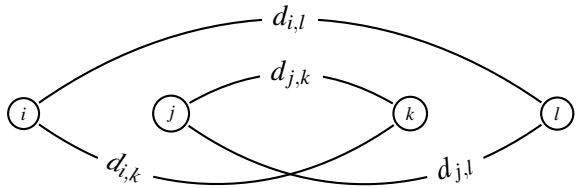
\begin{figure}[ht] 
\centering
\begin{tikzpicture}[thick,scale=0.6,-,shorten >= 1pt,shorten <= 1pt,every node/.style={scale=0.6}]
\begin{scope}[start chain=going right,node distance=15mm]
 \node[on chain,draw,circle] (m)  {$i$};
 \node[on chain,draw, circle] (k)  {$j$};
 \node[xshift=2cm,on chain,draw, circle] (j)  {$k$};
\node[on chain,draw,circle] (n) {$l$};
\end{scope}
\path[every node/.style={font=\sffamily\small}] (m) edge [bend left = 35]   node [fill=white] {$d_{i,l}$}  (n) ;
\path[every node/.style={font=\sffamily\small}] (k) edge [bend left]   node  [fill=white]  {$d_{j,k}$} (j) ;
\path[every node/.style={font=\sffamily\small}] (k) edge [bend right = 35]  node  [near end,sloped, fill=white] {$d_{j,l}$} (n) ;
\path[every node/.style={font=\sffamily\small}] (m) edge [bend right = 35]   node  [near start,sloped,fill=white] {$d_{i,k}$} (j) ;

\end{tikzpicture}
\vspace{-1ex}

\caption{\label{fig:4points}
 The convex quadrangle inequality.\label{4-point}}
\end{figure}
\vspace{-1ex}

\noindent 
A remarkable consequence of $(Q)$ is the existence of optimal transports that are {\em nested} in the sense that 
the flows do not intersect. This yields an efficient greedy algorithm to compute $d_{m,n}$. 
This algorithm is used in \S\ref{Sexamples} for various examples, revealing some particularly intriguing 
structures in the $d_{m,n}$'s.
We mention that the quadrangle inequality is closely related to the {\em inverse 
Monge property} \cite{bkr}, and
arises in various contexts with relevant algorithmic implications (see \cite{abkks,bei,bglz,knu,yao}).

\section{The recursive optimal transport metric}\label{smetric}

Let us start by proving that the $d_{m,n}$'s define a distance. This extends a result by \cite[Aygen-Satik]{aig} 
proved for the special case of the Krasnosel'skii iteration with  $\pi^{n}=(1-\alpha_{n})\pi^{n-1}+\alpha_{n}\delta^{n}$. 
The original proof was highly technical, 
though it was greatly simplified in \cite{brc} for the special case where $\alpha_n \geq \frac{1}{2}$ for all $n\in \NN$. 
This simpler proof is reproduced below, suitably adapted to the more general class of $\pi^n$'s considered here.
\begin{theorem}\label{t1}
The map $(m,n)\mapsto d_{m,n}$ is a metric on $\NNN=\NN\cup\{-1\}$.
\end{theorem}
\begin{proof}  
We already observed that $d_{m,n}=d_{n,m}$ and $d_{n,n}=0$, 
so we only need to establish the triangle inequality  and that $d_{m,n}>0$ for $m\neq n$.
We show inductively that these properties hold for $m,n\leq\ell$ for each $\ell\in \NNN$. 
The base case $\ell=-1$ is trivial. Suppose that both properties hold up to $\ell-1$ and let us prove them for $\ell$.

\vspace{1ex}

\noindent\underline{\em For $m,n\leq\ell$ with $m\neq n$ we have $d_{m,n}>0$}: 
Let $z$ be  an optimal transport for $d_{m,n}$. Since $\pi^m\neq\pi^n$ we can find $i\neq j$ with $z_{i,j}>0$
and the induction hypothesis implies $d_{m,n}\geq z_{i,j}\,d_{i-1,j-1}>0$. 

\vspace{1ex}

\noindent\underline{\em For $m,n,p\leq\ell$ we have $d_{m,n}\leq d_{m,p}+d_{p,n}$}: Let $z^{m,p}$ and $z^{p,n}$ be 
optimal transports for $d_{m,p}$ and $d_{p,n}$ respectively, and define
$$(\forall i=0,\ldots,m)(\forall j=0,\ldots,n)\quad z_{i,j}=\mbox{$\sum_{k=0}^p\omega_{i,k,j}$}$$
with 
$$\omega_{i,k,j}=\left\{\begin{array}{cl}
\dfrac{z^{m,p}_{i,k}z^{p,n}_{k,j}}{\pi_k^p}&\mbox{if }\pi_k^p\neq 0\\[1.5ex]
0&\mbox{otherwise.}
\end{array}\right.
$$
A straightforward computation shows that 
\begin{equation}\label{suma}
\left\{\begin{array}{rcl}
\mbox{$\sum_{i=0}^m\omega_{i,k,j}$}&=&z_{k,j}^{p,n}\\[1.5ex]
\mbox{$\sum_{j=0}^n\omega_{i,k,j}$}&=&z_{i,k}^{m,p}
\end{array}\right.
\end{equation}
from which it readily follows  that $z\in \Fcal_{m,n}$. Using the induction hypothesis 
and \eqref{suma} we get
$$\begin{array}{ccl}
d_{m,n}&\leq& \sum_{i=0}^m\sum_{j=0}^n z_{i,j}d_{i-1,j-1}\\[1.5ex]
&=& \sum_{i=0}^m\sum_{j=0}^n\sum_{k=0}^p \omega_{i,j,k}d_{i-1,j-1}\\[1.5ex]
&\leq& \sum_{i=0}^m\sum_{j=0}^n\sum_{k=0}^p \omega_{i,j,k}(d_{i-1,k-1}+d_{k-1,j-1})\\[1.5ex]
&=& \sum_{i=0}^m\sum_{k=0}^p z^{m,p}_{i,k}d_{i-1,k-1}+\sum_{j=0}^n\sum_{k=0}^p  z^{p,n}_{k,j}d_{k-1,j-1}\\[1.5ex]
&=& d_{m,p}+d_{p,n}.
\end{array}
$$

\vspace{-3ex}

\ \hfill\end{proof}

Now, part of the demand $\pi_i^n$ can be fulfilled at zero cost by shipping 
as much as possible  from the corresponding supply node $i$. This suggests 
the following type of simple transport plans. 
\begin{definition} A transport $z$ from $\pi^m$ to $\pi^n$ with $m\leq n$ is called {\em simple} if
$z_{i,i} = \min\{\pi_i^m,\pi_i^n\}$ for all $i=0,\ldots,m$.
\end{definition}
\begin{remark}\label{Remark1} If $z$ is a simple transport from $\pi^m$ to $\pi^n$, then when $\pi_i^m\geq \pi_i^n$ we have $z_{i,i}=\pi_i^n$ and
therefore $z_{j,i}=0$ for $j\neq i$. Similarly, if $\pi_i^m\leq \pi_i^n$ then $z_{i,i}=\pi_i^m$ and $z_{i,k}=0$ for $k\neq i$.
\end{remark}

\vspace{1ex}
As a consequence of the triangle inequality, it turns out that we may always restrict ourselves to transports that are simple.
\begin{proposition}\label{t2} 
Each $d_{m,n}$ admits a simple optimal transport. 
\end{proposition}
\begin{proof} Let $z$ be an optimal transport. If $z_{i,i}<\min\{\pi_i^m,\pi_i^n\}$ we must have $z_{i,k}>0$ for some $k\neq i$
and $z_{j,i}>0$ for some $j\neq i$. Decreasing  $z_{j,i}$ and $z_{i,k}$ by $\varepsilon$ while
increasing $z_{i,i}$ and $z_{j,k}$ by the same amount (see Figure \ref{figr}), 
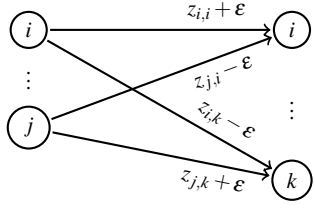
\begin{figure}[t]
\centering
\begin{tikzpicture}[thick,-,shorten >= 1pt,shorten <= 1pt,scale=0.5,every node/.style={scale=0.7}]
 \node[draw, circle] (i) at (0,-2) {$i$};
 \node[] (ii) at (0,-3.2) {$\vdots$};
 \node[draw, circle] (j) at (0,-4.6) {$j$};
 \node[draw,circle] (ii) at (7,-2) {$i$};
 \node[] (jj) at (7,-4) {$\vdots$};
 \node[draw, circle] (kk) at (7,-6) {$k$};
\draw[->] (i) -- (ii) node [above,near end] { $z_{i,i} +\varepsilon$};
\draw[->] (i) -- (kk) node [above, near end, sloped] { $z_{i,k} -\varepsilon$};
\draw[->] (j) -- (ii) node [below, near end, sloped] { $z_{j,i} -\varepsilon$} ;
\draw[->] (j) -- (kk) node [below, near end, sloped] { $z_{j,k} +\varepsilon$};
\end{tikzpicture}
\caption{\label{figr}Redistribution of flow for simple optimal transports}
\end{figure}
the modified transport is still feasible 
and the cost is reduced by
$$[d_{i-1,i-1}+d_{j-1,k-1}-d_{j-1,i-1}-d_{i-1,k-1}]\varepsilon\leq 0$$
so it remains optimal. Thus we can increase each $z_{i,i}$
up to $\min\{\pi_i^m,\pi_i^n\}$.
\end{proof}

\begin{figure}[t]
\centering
\begin{tikzpicture}[scale=0.65]
    \begin{axis}[
       ybar stacked,
        ymin=-0.2,
        ymax=0.18,
        xmin=-1.5,
        xmax=26.5,
        width  = 12cm,
        height = 8cm,
        bar width=9pt,
        xtick = data,
       xticklabels from table={ 0\\ \  \\ \  \\ \ \\ \ \\ 5\\ \ \\ \ \\ \ \\ \ \\ 10\\ \ \\ \ \\ \ \\ \ \\ 15\\ \ \\ \ \\ \ \\ \ \\20\\ \ \\ \ \\ \ \\ \ \\25\\}{[index]0},
       table/header=false,
        table/row sep=\\,
        enlarge y limits={value=0.2,upper},
    ]
 \addplot[fill=gray] table[x expr=\coordindex,y index=0]{ 
    0\\
    0\\
    0\\
    0\\
    0\\
    0\\
    0\\
    0\\
    0\\
    0\\
    0\\
    0\\
    0\\
    0\\
    0\\
    0\\
    0\\
    0\\
    0\\
    0\\
    0\\ 
    0.0760\\
    0.1251\\
    0.1397\\
    0.1826\\
    0.2360\\
    };
    \addplot[fill=black] table[x expr=\coordindex,y index=0]{ 
   -0.0022\\
   -0.0008\\
   -0.0008\\
   -0.0016\\
   -0.0018\\
   -0.0023\\
   -0.0040\\
   -0.0049\\
   -0.0059\\
   -0.0096\\
   -0.0115\\
   -0.0155\\
   -0.0223\\
   -0.0258\\
   -0.0314\\
   -0.0487\\
   -0.0759\\
   -0.0838\\
   -0.0933\\
   -0.1429\\
   -0.1747\\
    0\\
    0\\
    0\\
    0\\
    0\\
    };
      \end{axis}
\end{tikzpicture}
\caption{\label{pis2} Residual demands $\delta^{mn}=\pi^n-\pi^m$ ($m=20, n=25$).}
\end{figure}
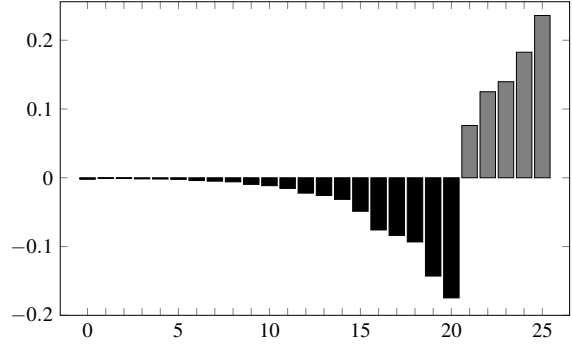

A further consequence of the fact that the transport cost $(i,j) \mapsto d_{i-1,j-1}$ is a distance, 
is the following Kantorovich-Rubinstein dual 
characterization for all $m\leq n$ (see {\em e.g.} \cite[Villani, Section 1.2]{vil}) where $\delta^{mn}_j=\pi^n_j - \pi^m_j$ denote the residual demands (see Figure \ref{pis2})
\begin{equation*}
d_{m,n}=\max_{u\in\RR^{n+1}}\left\{\mbox{$\sum_{j=0}^{n}$}\,\delta_j^{mn}\,u_j:u_j\leq u_i+ d_{i-1,j-1}\mbox{ $\forall i,j=0,\ldots,n$}\right\}. \leqno (\mathcal{D}_{m,n})
\end{equation*}
This is a refinement of linear programming duality, and each pair of
primal-dual optimal solutions $z^{mn}$ and $u^{mn}$ satisfy the complementary slackness 
 \begin{equation}\label{cs}
z_{i,j}^{mn}(u_j^{mn}\!-u_i^{mn})=z^{mn}_{ij}d_{i-1,j-1}\mbox{ for all $i,j=0,\ldots,n$}.
\end{equation}

\vspace{1ex}
\begin{remark}\label{Remark2} 
Every feasible $(u_i)^n_{i=0}$ in $({\mathcal D}_{m,n})$
can be extended  by setting $u_{i}={\displaystyle \min_{0\leq k\leq n}}u_{k}+d_{k-1,i-1}$ for $i>n$, 
so that the triangle inequality yields
\begin{equation}\label{modulo}
|u_{i}-u_{j}|\leq d_{i-1,j-1}\mbox{ for all }i,j\in\NN.
\end{equation}
This is a special case of the MacShane-Whitney extension of Lipschitz functions.
In particular all the $u_i$'s are within a distance at most 1 and,
since the objective function in $({\mathcal D}_{m,n})$ is invariant under 
translation, we may further assume that $u_i\in[0,1]$ for all $i\in\NN$. 
\end{remark}

\vspace{1ex}
\begin{remark}\label{Remark3} Under {\sc (h)} there is a unique simple
optimal transport from $\pi^n$ to $\pi^{n+1}$ with $z_{i,i}\!=\pi_i^{n+1}$ 
and $z_{i,n+1}\!=\pi_i^n\!-\!\pi_i^{n+1}\!$ for $0\leq i\leq n$. Moreover, 
setting $u_i=1-d_{i-1,n}\in[0,1]$ for all $0\leq i\leq n+1$, the triangle inequality implies that $u$
is a feasible solution for $({\mathcal D}_{n,n+1})$ which is also optimal because
$$\mbox{$\sum_{i=0}^{n+1}\delta_i^{n,n+1}u_i=\sum_{i=0}^nz_{i,n+1}d_{i-1,n}=d_{n,n+1}$.}$$
\end{remark}

\section{Tightness of the optimal transport bounds}\label{sfin}
Exploiting the dual characterization of the distances $d_{m,n}$ we 
proceed to establish that the optimal transport bounds in Theorem \ref{prop0} 
are the best possible estimates for the Krasnosel'skii-Mann iterates. 
Our results build upon similar ideas in \cite{bb92,brc}, suitably adapted to 
the more general $\pi^n$'s. 


\begin{theorem}\label{teoint}

Let $\D$ be the set of all  pairs of integers $(m,n)$ with $0\!\leq\! m\!\leq\! n$, and consider the unit cube $C=[0,1]^{\D}$ in the space $(\ell^\infty(\D),\|\cdot\|_\infty)$.
Then, for each sequence $(\pi^n)_{n\in\NN}$ there exists  a non-expansive map $T:C\to C$ 
and a corresponding {\sc (km)} sequence with $\|x^m-x^n\|_\infty= d_{m,n}$ for all $m,n\in\NN$.
Moreover, under {\sc (h)} we also have the tight bound $\|x^n-Tx^n\|_\infty=R_n$.
\end{theorem}

\begin{proof}
For each $(m,n)\in\D$ consider an optimal solution $u^{mn}$  for $({\mathcal D}_{m,n})$ and 
 its extension as in Remark \ref{Remark2} so that $u^{mn}_i\in [0,1]$ and \eqref{modulo} hold for all $i,j\in\NN$.

For every integer $k\in\NN$ define $y^{k}\in C=[0,1]^\D$ as
 \begin{equation}
\label{yyy}  y^k_{m,n}=u^{mn}_{k}\qquad\forall (m,n)\in \D,
\end{equation}
and  a corresponding sequence $x^{k}\in C$ for $k\geq 0$ given by
 \begin{equation}
 \label{rec} \mbox{$ x^k=\sum_{i=0}^k\pi_i^ky^{i}$}
 \end{equation}
 so that in particular $x^0=y^{0}=(u_0^{mn})_{(m,n)\in \D}$.
 
 We claim that $\|y^{m+1}-y^{n+1}\|_\infty\leq d_{m,n} =\|x^m-x^n\|_\infty$ for all $m,n\in\NN$.
 Indeed, from \eqref{modulo} we get $|u^{m'n'}_{m+1}\!-u^{m'n'}_{n+1}|\leq d_{m,n}$ for all $(m',n')\in \D$, 
which implies $$\|y^{m+1}\!-y^{n+1}\|_\infty=\|(u^{m'n'}_{m+1}\!-u^{m'n'}_{n+1})_{(m',n')\in \D}\|_\infty\leq d_{m,n}.$$
Also, selecting an optimal transport $z^{mn}$ for $({\mathcal P}_{m,n})$ and proceeding as in the proof of Theorem \ref{prop0} 
we have
\begin{equation}\label{diferencia}
\mbox{$x^m-x^n=\sum_{i=0}^{m}\sum_{j=0}^{n}z^{mn}_{i,j}(y^{i}-y^{j})$}
\end{equation}
so that the triangle inequality yields
\begin{equation}\label{siete}
\mbox{$
\|x^m\!-x^n\|_\infty\leq 
\sum_{i=0}^{m}\sum_{j=0}^{n}z^{mn}_{i,j}d_{i-1,j-1}=d_{m,n}
$}
\end{equation}
while considering the $(m,n)$-coordinate in \eqref{diferencia}, the complementary slackness 
\eqref{cs} gives
 \begin{eqnarray*}
 |x_{m,n}^m-x_{m,n}^n|&=&|\mbox{$\sum_{i=0}^{m}\sum_{j=0}^{n}z^{mn}_{i,j}(y_{m,n}^{i}-y_{m,n}^{j})$}|\\
 &=&| \mbox{$\sum_{i=0}^{m}\sum_{j=0}^{n}z^{mn}_{i,j}(u^{mn}_i-u^{mn}_j)$}|\\
 &=&\mbox{$\sum_{i=0}^{m}\sum_{j=0}^{n}z^{mn}_{i,j}d_{i-1,j-1}=d_{m,n}$}
 \end{eqnarray*} 
 which combined with \eqref{siete} yields $\|x^m-x^n\|_\infty=d_{m,n}$ as claimed. 
 
Define $T:S\to C$ on the set $S=\{x^k:k\in\NN\}\subseteq C$ by $Tx^k=y^{k+1}$, so that $T$
is non-expansive. Since $\ell^\infty(\D)$ as well as the unit cube $C$ are hyperconvex, by Theorem 4 in
 \cite[Aronszajn \& Panitchpakdi]{aro}, $T$ can be extended to a non-expansive map $T:C\to C$ and then
\eqref{rec} is precisely a {\sc (km)} sequence which attains all the bounds 
$\|x^m-x^n\|_\infty=d_{m,n}$ with equality.

It remains to prove that under {\sc (h)} we have $\|x^n-Tx^n\|_\infty=R_n$.
The upper bound follows from \eqref{residual} or by direct computation
$$\|x^n-Tx^n\|_\infty=\mbox{$\|\sum_{i=0}^n\pi^n_i(y^{i}-y^{n+1})\|_\infty$}\leq \mbox{$\sum_{i=0}^n\pi_i^nd_{i-1,n}=R_n$}.$$
For the reverse inequality, let $u^{n,n+1}$ be chosen as in Remark \ref{Remark3},
so that
$u^{n,n+1}_{n+1}-u_i^{n,n+1}=d_{i-1,n}$ for each $0\leq i\leq n$, and therefore 
\begin{eqnarray*}
\|x^n-Tx^n\|_\infty&=&\mbox{$\|\sum_{i=0}^n\pi^n_i(y^{i}-y^{n+1})\|_\infty$}\\
&\geq&|\mbox{$\sum_{i=0}^n\pi^n_i(y_{n,n+1}^{i}-y_{n,n+1}^{n+1})|$}\\
&=&|\mbox{$\sum_{i=0}^n\pi^n_i(u_i^{n,n+1}-u^{n,n+1}_{n+1})|$}\\
&=&\mbox{$\sum_{i=0}^n\pi^n_id_{i-1,n}=R_n$}
\end{eqnarray*}
which completes the proof.
 \end{proof}
 
 At this point we do not know if $R_n$ is still a sharp bound for the residuals when {\sc (h)} fails,
or whether a tighter bound might hold in that case.
On the other hand, we observe that the unit cube $C=[0,1]^\D$ is not only hyperconvex but also bounded
so that, from a result by \cite[Sine]{sin} and \cite[Soardi]{soa},
the map $T$ in Theorem \ref{teoint} has fixed points, that is, $\mathop{\rm Fix}(T)$
is nonempty. For a survey of the extremely rich structure and results on hyperconvex spaces we refer to 
\cite[Esp\'inola \& Khamsi]{esp}.

\section{Krasnosel'skii and Halpern iterations}\label{s51}

Let us illustrate the previous tight estimates for the classical Krasnosel'skii and Halpern iterations.
For the Krasnosel'skii iteration $$x^{n}=(1-\alpha_n)x^{n-1}+\alpha_n Tx^{n-1},\qquad \alpha_n\in(0,1)$$
we have $\pi^n=(1-\alpha_n)\pi^{n-1}+\alpha _n\delta^n$ and {\sc (h)} holds automatically. In view of 
Remark \ref{Remark3} the unique simple optimal transport from $\pi^n$ to $\pi^{n+1}$ gives
$$\mbox{$d_{n,n+1}=\sum_{i=0}^n(\pi_i^n-\pi_i^{n+1})d_{i-1,n}=\alpha_{n+1}\sum_{i=0}^n\pi_i^n d_{i-1,n}=\alpha_{n+1}R_n$}$$
which yields a simpler expression for the fixed point residual bound
 $$\|x^n-Tx^n\|\leq R_n=\frac{d_{n,n+1}}{\alpha_{n+1}}.$$ 
 
 This coincides with the tight bound established in \cite[Bravo \& Cominetti]{brc}. 
 Moreover, in view of \cite[Cominetti {\em et al.}]{csv}
 (see also \cite[Baillon \& Bruck]{bb96} for the case when $\alpha_k\equiv\alpha$),  we have
 \begin{equation}\label{bound}
 R_n=\frac{d_{n,n+1}}{\alpha_{n+1}}\leq\frac{\mathop{\rm diam}(C)}{\sqrt{\pi\sum_{k=1}^n\alpha_k(1\!-\!\alpha_k)}}
 \end{equation}
 which converges to 0 as soon as $\sum_{k=1}^\infty\alpha_k(1\!-\!\alpha_k)=\infty$. In particular, for 
 $\alpha_k\equiv\alpha$ we recover the rate $R_n\sim O(1/\sqrt{n})$ proved in \cite{bb96}. 
 
 Let us mention that for $\alpha_k\equiv\alpha\approx 1$ the bound \eqref{bound} becomes
 asymptotically tight   and the proportionality constant $1/\sqrt{\pi}$ is the best possible (see  \cite{brc}).
Recall also that when $X$ is uniformly convex  with Fr\'echet differentiable norm,
 the property $\sum_{k=1}^\infty\alpha_k(1\!-\!\alpha_k)=\infty$ implies weak convergence of the iterates 
 $x^n$ to a fixed point, provided $\mathop{\rm Fix}(T)$ is nonempty (see \cite[Reich]{rei}).
 
 \vspace{1ex}
 Let us consider next Halpern's iteration 
$$x^n=(1-\beta_n)y^0+\beta_nTx^{n-1},\qquad \beta_n\in(0,1).$$
In this case  we have the simpler structure $\pi_0^n=1-\beta_n$ and $\pi_n^n=\beta_n$ so that
 \begin{equation}\label{aux01}
 \mbox{$R_n=\sum_{i=0}^n\pi_i^nd_{i-1,n}=(1-\beta_n)+\beta_nd_{n-1,n}.$}
 \end{equation}
Here {\sc (h)} amounts to $\beta_n$ being non-decreasing, in which case 
 Remark \ref{Remark3} yields
 $$d_{n-1,n}=(\beta_n-\beta_{n-1})+\beta_{n-1} d_{n-2,n-1}.$$
 This combined with \eqref{aux01} implies $R_n=(1-\beta_n)^2+\beta_n R_{n-1}$ which has an
explicit solution: letting $B_i^n=\prod_{k=i}^n\beta_k$ with $\beta_0=0$ and $B_{n+1}^n=1$,
we get 
 $$R_n=\mbox{$\sum_{i=0}^n(1-\beta_i)^2B_{i+1}^n.$}$$
 
For the particular case $\beta_n=\frac{n}{n+1}$ a direct calculation gives
\begin{equation}\label{halrate}
\mbox{$\|x^n-Tx^n\|\leq R_n=\frac{H_{n+1}}{n+1}\sim O(\frac{\ln n}{n})$}
\end{equation}
with $H_n=\sum_{k=1}^n\frac{1}{k}$ the $n$-th harmonic number. From Theorem \ref{teoint} we know that this bound 
is sharp and cannot be improved. However,  \cite[Lieder]{lie}   recently proved that when $(X,|\cdot|)$ is a Hilbert space and 
$T:X\to X$ has a fixed point $x^*=Tx^*$, 
the choice $\beta_n=\frac{n}{n+1}$ yields the faster rate $$\mbox{$|x^n-Tx^n|\leq \dfrac{2|x^0-x^*|}{n+1}\sim O(\frac{1}{n}).$}$$
In order to compare with \eqref{halrate}, consider the rescaled norm $\|\cdot\|=\frac{1}{r_0}|\cdot|$ with $r_0=2|x^0-x^*|$ 
so that $C=B(x^*,|x^0-x^*|)$ has diameter $1$. Then
Lieder's bound reads $\|x^n-Tx^n\|\leq \frac{1}{n+1}$ which is clearly smaller than \eqref{halrate}. 
This does not contradict Theorem \ref{teoint}: it simply shows that the tight bounds are attained in
a non-Hilbert setting such as the space $\ell^\infty(\D)$.

On the other hand, a minor modification of Halpern with $\beta_n=\frac{n}{n+2}$ yields the faster rate
$R_n=\frac{4}{n+1}(1-\frac{H_{n+2}}{n+2})\sim O(\frac{1}{n})$ in general normed spaces. 
This  bound is also tight and   slightly  improves the estimate $\|x^n-Tx^n\|\leq \frac{4}{n+1}$ established in 
\cite[Sabach \& Shtern, Lemma 5]{sab}. For the record, we observe that although $\beta_n=\frac{n}{n+2}$ achieves a faster 
rate than $\beta_n=\frac{n}{n+1}$, the slow growth of $H_n$ implies that for small $n$ both sharp bounds coincide 
within a moderate  multiplicative factor (less than 2.5 up to $n=10.000$).

In general it remains open to find conditions under which $R_n\to 0$.
The forthcoming paper \cite[Contreras \& Cominetti]{cc} further investigates various algorithms 
that fit the general framework {\sc (km)}, and determines the optimal choices for the 
weights $\pi^n$ that guarantee the faster rate $R_n\sim O(1/n)$ with the best  
proportionality constants.

\section{Monotonicity, convex quadrangle inequality, and the inside-out greedy algorithm}
\label{squadrangle}

In this section we discuss some further properties of the optimal transport metric 
under the additional assumption {\sc (h)}. A first 
consequence is that the distances $d_{m,n}$ increase as 
$m$ and $n$ get farther apart: for any fixed $m$ the map $n\mapsto d_{m,n}$ 
decreases from a value of 1 at $n=-1$ down to 0 at $n=m$, after which it increases 
for $n\geq m$. 
Intuitively, this comes from the fact that the 
distribution $\pi^n$  drifts gradually towards larger integers, so that for $n\geq m$ 
it is more expensive to transport $\pi^m$ towards $\pi^{n+1}$ than it is to 
satisfy the demands $\pi^n$ which are  located somewhat closer to $\pi^m$. 
A similar argument can be made for $n<m$.
We formally prove this monotonicity by exploiting 
both the primal and dual caracterizations of $d_{m,n}$.

\begin{theorem}\label{t3}
Assume {\sc (h)}. Then in the region $m\leq n$ we have that $d_{m,n}$ decreases with $m$ and increases with $n$.
\end{theorem}
\begin{proof} Let us first show that $d_{m,n}\leq d_{m-1,n}$ for $0\leq m\leq n$. Since for $n=0$ this holds trivially,
we proceed by induction assuming that these inequalities hold up to $n-1$. To establish the property for $n$ we make a second
induction on $m$.
For $m=0$ we clearly have $d_{0,n}\leq d_{-1,n}=1$. Now, take $z$ a simple optimal transport 
 for $d_{m-1,n}$ so that for $i=0,\ldots,m-1$ we have $z_{i,i}=\pi_i^n$ with the excess supply $\pi^{m-1}_i\!-\pi^n_i$ 
 shipped to nodes $j\geq m$.
 We transform $z$ into a feasible transport for $d_{m,n}$ while reducing the cost.
 To this end, for each $i=0,\ldots,m-1$ with $\pi_{i}^{m}<\pi_i^{m-1}$ we take some positive flow $z_{i,j}>0$ with 
 $m\leq j\leq n$, from which we remove $\varepsilon>0$ and increase $z_{m,j}$ by the same amount (see Figure \ref{figflow2}).
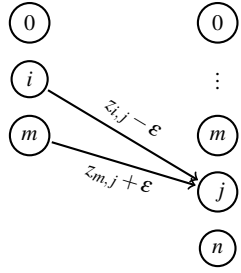
\begin{figure}[t]
\centering
\begin{tikzpicture}[thick,-,shorten >= 1pt,shorten <= 1pt,scale=0.5,every node/.style={scale=0.7}]
 \node[draw,circle] (ii0) at (5,0) {$0$};
 \node[] (idost)  at (5,-1.4)  {$\vdots$};
\node[draw,circle] (iin1) at (5,-3)  {$m$};
 \node[draw, circle,text width=1mm] (ii1) at (5,-4.5) {$j$};
\node[draw,circle] (iim1) at (5,-6) {$n$};

 \node[draw,circle] (jj0) at (0,0)  {$0$};
 \node[draw, circle] (jjaux) at (0,-1.5) {$i$};
 \node[draw,circle] (jj2) at (0,-3) {$m$};

\draw[<-] (ii1) -- (jjaux)  node [above, midway,sloped]{ $z_{i,j} -\varepsilon$};
\draw[<-] (ii1) -- (jj2) node [below, midway, sloped] { $z_{m,j} +\varepsilon$} ;
\end{tikzpicture}
\caption{\label{figflow2} Redistribution of flows to prove monotonicity.} 
\end{figure}
Since $i-1< m-1\leq j-1\leq n-1$, the induction hypothesis implies that the cost decreases by
$[d_{m-1,j-1}-d_{i-1,j-1}]\varepsilon\leq 0$. We repeat these flow transfers until the outflow $\sum_{j=0}^nz_{i,j}$ 
 is reduced to $\pi_i^m$.  At this point a fraction $\pi_i^{m-1}\!-\pi_i^m$ of the excess suply of each $i\leq m-1$ 
 has been fully transferred to $m$ and therefore 
 $$\mbox{$\sum_{j=0}^nz_{m,j}=\sum_{i=0}^{m-1} (\pi_i^{m-1}\!-\pi_i^m)=1- \sum_{i=0}^{m-1} \pi_i^m=\pi_m^m.$}$$ 
Hence, the final flow is feasible for $d_{m,n}$ and since the cost was decreased at every step 
we conclude $d_{m,n}\leq C_{m,n}(z)\leq d_{m-1,n}$.

We now use the dual characterization $({\mathcal D}_{m,n})$ to prove $d_{m,n}\leq d_{m,n+1}$ by induction on $n$. This holds trivially for $n=m$.
Assume that it holds up to $n-1$ and let us prove it for $n$.
Take $u^{mn}$ optimal for  $({\mathcal D}_{m,n})$ and consider the vector $u$ with 
$u_i=u^{mn}_i$ for $i=0,\ldots,n$ and $u_{n+1}=u^{mn}_{j_0}$ where
$u^{mn}_{j_0}$ is the maximal coordinate in $u^{mn}$. By construction 
$|u_i-u_j|\leq d_{i-1,j-1}$ for all $0\leq i,j\leq n$ while the induction hypothesis gives
$$|u_i-u_{n+1}|=|u^{mn}_i-u^{mn}_{j_0}|\leq d_{i-1,j_0-1}\leq d_{i-1,n}$$
so that  $u$ is feasible in $({\mathcal D}_{m,n+1})$. 
On the other hand, 
$\delta^{m,n}_j=\delta^{m,n+1}_j+\varepsilon_j$ with $\varepsilon_j=\pi^{n}_j-\pi^{n+1}_j\geq 0$  for $j=0,\ldots, n$ 
from which we get
$$\mbox{$d_{m,n}=\sum_{j=0}^n\delta^{m,n}_ju_j=\sum_{j=0}^n\delta^{m,n+1}_ju_j+\sum_{j=0}^n\varepsilon_ju_j$}.$$
Since $u_j\leq u_{n+1}$ this last sum can be majorized by $\sum_{j=0}^n\varepsilon_ju_{n+1}$, which combined with
the equality $\sum_{j=0}^n\varepsilon_j=\pi_{n+1}^{n+1}=\delta^{m,n+1}_{n+1}$ implies
$$\mbox{$d_{m,n}\leq\sum_{j=0}^{n+1}\delta^{m,n+1}_ju_j\leq d_{m,n+1}$}.$$
completing the induction step.
 \end{proof}

\vspace{1ex}
\noindent
{\sc Remark 4}. The informal argument just before Theorem \ref{t3} might suggest that monotonicity could 
hold under the weaker condition of first order stochastic dominance 
$\pi^n\preceq_{st}\pi^{n+1}$, that is $\sum_{i\geq k}\pi_i^n\leq \sum_{i\geq k}\pi_i^{n+1}$ for all $k\in\NN$.
It turns out that this is not always the case as shown by 
\begin{eqnarray*}
\pi^0&=&\mbox{$(1,0,0,0,\ldots )$}\\
\pi^1&=&\mbox{$(\frac{4}{5},\frac{1}{5},0,0,\ldots )$}\\
\pi^2&=&\mbox{$(\frac{1}{2},\frac{1}{4},\frac{1}{4},0,0,\ldots )$}\\
\pi^3&=&\mbox{$(\frac{1}{2},\frac{1}{4},0,\frac{1}{4},0,0,\ldots )$}\\
\pi^4&=&\mbox{$(0,\frac{1}{2},\frac{1}{4},0,\frac{1}{4},0,0,\ldots )$}\\
\pi^5&=&\mbox{$(0,0,\frac{1}{2},0,0,\frac{1}{2},0,0,\ldots )$}
\end{eqnarray*}
for which we have $d_{3,5}=\frac{22}{40}$ and $d_{4,5}=\frac{23}{40}$ so that $d_{3,5}<d_{4,5}$.

\vspace{2ex}
We next show that under {\sc (h)} the distances $d_{m,n}$ satisfy  the 
convex quadrangle inequality $(Q)$. The proof proceeds inductively, by considering 
the quadrangle inequality up to $n$, namely
\begin{equation*}\label{4pointn}
d_{i,l}+d_{j,k}\leq d_{i,k}+d_{j,l}\qquad   \mbox{ for all $-1\leq i\leq j \leq k\leq l\leq n$.}
\leqno (Q_{n})
\end{equation*}
A key consequence of the latter inequalities is the existence 
of optimal transports that are not only simple but also nested in the sense that the flows do not intersect.
\begin{definition} A transport plan $z$ is called {\em nested}
if there are no $i<j<k<l$ with $z_{i,k}>0$ and $z_{j,l}>0$.
\end{definition}

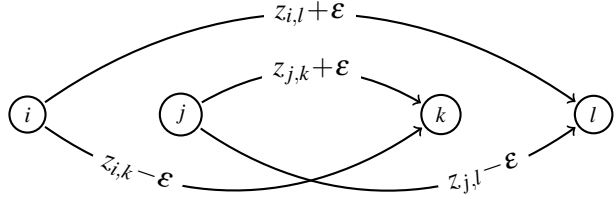
\begin{figure}[t] 
\centering
\begin{tikzpicture}[thick,-,shorten >= 1pt,shorten <= 1pt,scale=0.3,every node/.style={scale=0.7}]
\begin{scope}[start chain=going right,node distance=15mm]
 \node[on chain,draw,circle] (m)  {$i$};
 \node[on chain,draw, circle] (k)  {$j$};
 \node[xshift=2cm,on chain,draw, circle] (j)  {$k$};
\node[on chain,draw,circle] (n) {$l$};
\end{scope}
\path[every node/.style={font=\sffamily\small}] (m) edge [->,bend left = 35]   node [fill=white] {$z_{i,l}\!+\!\varepsilon$}  (n) ;
\path[every node/.style={font=\sffamily\small}] (k) edge [->,bend left]   node  [fill=white]  {$z_{j,k}\!+\!\varepsilon$} (j) ;
\path[every node/.style={font=\sffamily\small}] (k) edge [->,bend right = 35]  node  [near end,sloped, fill=white] {$z_{j,l}\!-\!\varepsilon$} (n) ;
\path[every node/.style={font=\sffamily\small}] (m) edge [->,bend right = 35]   node  [near start,sloped,fill=white] {$z_{i,k}\!-\!\varepsilon$} (j) ;
\end{tikzpicture}
\caption{\label{figr2}Recirculation of $\varepsilon$ along the cycle $i\rightarrow l\rightarrow j \rightarrow k \rightarrow i$.}
\end{figure}

\begin{proposition}\label{Le1} Assume $(Q_{n-1})$.
Then each $d_{m,n}$ for $m\leq n$ admits an optimal transport that is both simple and nested. 
\end{proposition}
\begin{proof}
Let $z$ be a simple optimal transport for $d_{m,n}$ and suppose that $z_{i,k}>0$ and $z_{j,l}>0$ for some 
$0\leq i<j<k< l\leq n$. 
Removing $\varepsilon=\min\{z_{i,k},z_{j,l}\}$  from $z_{i,k}$ and $z_{j,l}$  
while increasing $z_{i,l}$ and $z_{j,k}$ by the same amount  (see Figure~\ref{figr2}),
the modified transport remains feasible and $(Q_{n-1})$ implies a cost reduction of
$$[d_{i-1,l-1}+d_{j-1,k-1}-d_{i-1,k-1}-d_{j-1,l-1}]\varepsilon\leq 0.$$
Hence, the modified transport remains optimal with $z_{i,k}=0$ or $z_{j,l}=0$. 
Inductively we can eliminate all flow intersections and obtain an optimal transport that is simple and nested.
\end{proof}

Under {\sc (h)} the properties of being simple and nested determine a {\em unique transport plan} through the 
following greedy procedure. Imagine the supply nodes as a set of buckets $S_0,\ldots,S_m$ initially filled with 
volumes $\pi_i^m$ for $0\leq i\leq m$, and a set of empty  demand buckets $D_0,\ldots,D_n$ with capacities 
$\pi_j^n$ for $0\leq j\leq n$. Considering the supplies  $S_m,S_{m-1}\ldots,S_0$ in reverse order, we use the 
volume $\pi_i^m$ in each $S_i$ to fill $D_i$ and distribute the surplus $\sigma_i=\pi_i^m-\pi_i^n$ among the 
unmet demands $D_k$ for  $k> m$ in increasing order. Thus, $S_m$ fills $D_m,D_{m+1},\ldots, D_{k_m}$ 
where $D_{k_m}$ is the bucket at which the volume $\pi_m^m$ is depleted. Next $S_{m-1}$ fills $D_{m-1}$ 
and resumes at the bucket $D_{k_m}$  partially filled in the previous round and up to some $D_{k_{m-1}}$ at 
which $\pi_{m-1}^m$ is depleted, and so on.

%
%
%

More formally, the procedure is as follows:
{\small 
\begin{algorithm}[H]
\caption{ Inside-Out --- input $\pi^m,\pi^n$ with $m\leq n$}\label{alg:InsideOut}
\begin{algorithmic}[1]
\State For each $i=0,\ldots,m$ set $z_{i,i}=\pi_i^n$ and $z_{j,i}=0$ for $j\neq i$,
and compute the residual supply $\sigma_i=\pi_i^m-\pi_i^n$.\vspace{0.5ex}
\State Use $\sigma_m$ to fill the demands $\pi_{m+1}^n,\pi_{m+2}^n,\ldots$ 
up to the node $k_m$ at which $\sigma _m$ is depleted with $\pi_{k_m}^n$ only 
partially filled. This sets the flows $z_{m,j}=\pi_j^n$ for $j\in\{m\!+\!1,\ldots,k_m\!-\!1\}$ and $z_{m,k_m}=\sigma_m-\sum_{j=m+1}^{k_m-1}\pi_j^n$.\vspace{0.5ex}
\State Use $\sigma_{m-1}$ to fill any unmet demand at $k_m$ as well as the subsequent demands
up to a  node $k_{m-1}$ at which $\sigma _{m-1}$ is depleted. This determines the nonzero flows
$z_{m-1,j}$ for $j\in\{k_m,\ldots,k_{m-1}\}$.\vspace{0.5ex}
\State Proceed sequentially using $\sigma_i$ for $i=m\!-\!2,\ldots,0$
to fulfill the unmet demands at the nodes $j\in \{k_{i+1},\ldots,k_i\}$, which determine the corresponding  
 flows $z_{i,j}$. \vspace{0.5ex}
\end{algorithmic}
\end{algorithm}
}

As it turns out, the validity of  the quadrangle inequality and that of the inside-out algorithm are
closely related, and will be proved simultaneously. For the special case where 
$\pi_i^n= \alpha_i \prod_{k=i+1}^n(1-\alpha_k)$, the optimality of this algorithm was 
conjectured in \cite{bb92}, and later confirmed in \cite{aig} with a long and technical proof. 
Below we present a simpler argument that exploits the optimal transport structure of the $d_{m,n}$'s, 
and which applies to the more general distributions $\pi^n$'s considered here. 
The proof uses Proposition \ref{Le1} inductively  to show that $(Q_n)$ 
--- or more precisely, its equivalent form $(\tilde Q_n)$ below --- holds for all $n\in\NN$.
\begin{lemma}\label{Le3} Let $\Delta_{m,j}\triangleq d_{m,j+1}-d_{m,j}$. Then $(Q_n)$ is equivalent to 
$$ \mbox{For all $m<j<n$ we have $\Delta_{m,j}\leq \Delta_{m+1,j}$.}\leqno(\tilde Q_n)$$
\end{lemma}
\begin{proof} The inequality in $(Q_n)$ can be written as $d_{i,l}-d_{i,k}\leq d_{j,l}-d_{j,k}$
which amounts to the fact that $d_{m,l}-d_{m,k}$ increases with $m$ 
(for $m\leq k\leq l \leq n$). Since $d_{m,l}-d_{m,k}=\sum_{j=k}^{l-1}\Delta_{m,j}$
this is in turn equivalent to $(\tilde Q_n)$.
\end{proof}

\begin{theorem}\label{t4} Assume {\sc (h)}. Then, the recursive optimal transport metric $d_{m,n}$ satisfies the convex 
quadrangle inequality $(Q)$, and therefore the inside-out algorithm computes an optimal transport for all $m\leq n$.
\end{theorem}
\begin{proof} Since $(\tilde Q_0)$ holds trivially, it suffices to prove that $(\tilde Q_{n})\Rightarrow (\tilde Q_{n+1})$. 
In view of the induction hypothesis $(\tilde Q_{n})$, in order to prove $(\tilde Q_{n+1})$ it suffices 
to consider $j=n$ and to show that $\Delta_{m,n}\leq\Delta_{m+1,n}$ for all $m<n$.

Using $(\tilde Q_n)$, and according to Proposition \ref{Le1} and Lemma \ref{Le3}, it follows that the four 
problems $d_{m+1,n+1}$, $d_{m+1,n}$, $d_{m,n+1}$, $d_{m,n}$, admit 
simple and nested optimal transports given by the inside-out algorithm. 
Let $z$ and $z'$ be such optimal transports for $d_{m+1,n+1}$ and $d_{m+1,n}$ respectively. 
We will modify $z$ and $z'$ simultaneously, in such a way that $\Delta_{m+1,n}$ 
decreases and  the transports are progressively transformed into simple and nested 
transports for $d_{m,n+1}$ and $d_{m,n}$. By Proposition \ref{Le1} these modified transports 
are optimal for $d_{m,n+1}$ and $d_{m,n}$,
which yields the conclusion $\Delta_{m,n}\leq\Delta_{m+1,n}$.

The modification of $z$ and $z'$ is done in a sequence of stages.
The informal idea is as follows.
Since $\pi_{m+1}^{m+1}>\pi_{m+1}^m=0$, the initial transport $z$ has a surplus of 
flows $z_{m+1,j}$ shipped out from $m+1$, so that to regain feasibility for 
$d_{m,n+1}$ these must be substituted with flows $z_{i,j}$ from sources $i\leq m$. 
The same is needed to transform $z'$ into a feasible transport for $d_{m,n}$.
We do this in a series of flow shifts that progressively eliminate the outflow surplus 
at the tail supply nodes, pushing the imbalance towards sources closer  to 0
and iterating until all nodes dispatch the correct supplies $\pi_i^m$.
Along the process we make sure that the transports remain simple and
nested, while the quantity $\Delta_{m+1,n}$ is reduced at each step.

In order to formalize this idea, let $V_i=\sum_{k=i}^{m+1}\pi_k^{m+1}$ be the cumulative supply from $i$ to $m+1$.
Then the node $k_i>m+1$ at which $V_i$ is depleted in the inside-out procedure for 
$d_{m+1,n+1}$ is characterized by the inequalities
$$\sum_{j=i}^{k_i-1}\pi_j^{n+1}<V_i\leq \sum_{j=i}^{k_i}\pi_j^{n+1}.$$
Similarly, the corresponding node $k'_i>m+1$ for  $d_{m+1,n}$ is such that 
$$\sum_{j=i}^{k'_i-1}\pi_j^{n}<V_i\leq \sum_{j=i}^{k'_i}\pi_j^{n}.$$
From this it follows that $k'_i\leq k_i$. Indeed, this inequality is trivial if $k_i=n+1$, whereas
if $k_i\leq n$ we may use the inequalities $\pi_j^{n+1}\leq\pi_j^n$ to obtain
\begin{equation}\label{orderkm}
 \sum_{j=i}^{k'_i-1}\pi_j^n< V_i\leq \sum_{j=i}^{k_i}\pi_j^{n+1} \leq \sum_{j=i}^{k_i}\pi_j^{n}
 \end{equation}
from where we deduce that $k'_i-1<k_i$, and therefore $k'_i\leq k_i$ as claimed.

In a first phase of the transformation of $z$ we consider the supply $\pi_{m+1}^{m+1}$ that is shipped from $m+1$ 
towards the demand nodes $m+1,\ldots,k_{m+1}$. We transfer a small amount $\varepsilon$ from $m+1$ to $m$ (see Figure \ref{phase1}), 
which implies a change in the cost of $[d_{m-1,k_{m+1}-1}-d_{m,k_{m+1}-1}]\varepsilon$. An analog flow transfer in $z'$
induces a corresponding cost change of $[d_{m-1,k'_{m+1}-1}-d_{m,k'_{m+1}-1}]\varepsilon$. Both changes combined imply that 
$\Delta_{m+1,n}$ is decreased by  
$$[d_{m-1,k_{m+1}-1}-d_{m,k_{m+1}-1}-d_{m-1,k'_{m+1}-1}+d_{m,k'_{m+1}-1}]\varepsilon\leq 0$$
where negativity follows from $(\tilde Q_n)$ since $m-1\leq m\leq k_{m+1}'-1\leq k_{m+1}-1$. 
Let  $z_\varepsilon$, $z'_\varepsilon$ and 
$\Delta^\varepsilon_{m+1,n}=C_{m+1,n+1}(z_\varepsilon)-C_{m+1,n}(z'_\varepsilon)$ denote the modified transports and value. 

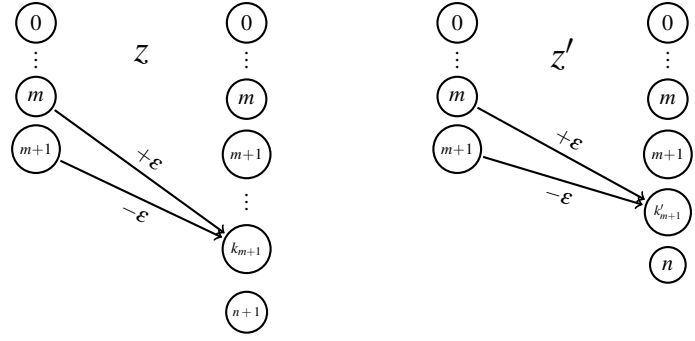
\begin{figure}[t]
\centering
\begin{tikzpicture}[thick,-,shorten >= 1pt,shorten <= 1pt,scale=0.7,every node/.style={scale=0.7}]
 \node[scale=1.4] (z) at (2,-1) {\large $z$};
 \node[scale=1.4] (zp) at (10,-1) {\large $z'$};
 \node[draw,circle] (i0) at (0,-0.4) {$0$};
 \node[] (i1) at (0,-1)   {$\vdots$};
\node[draw,circle] (im1) at (0,-1.8) {$m$};
\node[draw,circle, scale=0.7] (im11) at (0,-2.8) {$m\!+\!1$};
 \node[draw,circle] (j0) at (4,-0.4) {$0$};
 \node[] (j1) at (4,-1) {$\vdots$};
\node[draw,circle] (jm1) at (4,-1.85) {$m$};
 \node[draw, circle,scale=0.7] (jaux)  at (4,-2.9)  {$m\!+\!1$};
 \node[circle] (jauxxx)  at (4,-3.7)  {$\vdots$};
 \node[draw, circle,scale=0.7] (jaux1)  at (4,-4.7)  {$k_{m+1}$};
 \node[draw,circle, scale=0.6] (j2)  at (4,-5.9) {$n+1$};
\draw[->] (im1) -- (jaux1)  node [above, midway, sloped] { $+\varepsilon$} ;
\draw[->] (im11) -- (jaux1) node [below, midway,sloped] { $-\varepsilon$} ;
 \node[draw,circle] (ii0) at (8,-0.4)  {$0$};
 \node[] (ii1) at (8,-1)  {$\vdots$};
\node[draw,circle] (iim1) at (8,-1.8) {$m$};
\node[draw,circle, scale=0.7] (imm11) at (8,-2.8) {$m\!+\!1$};
 \node[draw,circle] (jj0)  at (12,-0.4) {$0$};
 \node[] (jj1)   at (12,-1) {$\vdots$};
\node[draw,circle] (jjm1)  at (12,-1.85) {$m$};
\node[draw,circle, scale=0.7] (jjaux)  at (12,-2.9) {$m\!+\!1$};
 \node[draw, circle,scale=0.7] (jjaux1)  at (12,-4.0) {\small $k'_{m+1}$};
 \node[draw,circle] (jj2)  at (12,-5.0) {$n$};
\draw[->] (iim1) -- (jjaux1)  node [above, midway,sloped]{ $+\varepsilon$};
\draw[->] (imm11) -- (jjaux1) node [below, midway, sloped] { $-\varepsilon$} ;
\end{tikzpicture}
\caption{\label{phase1} Phase 1 of redistribution of flows.} 
\end{figure}

As we proceed with these flow transfers the residual supply $\pi_{m+1}^{m+1}-\varepsilon$ of node 
$m+1$ decreases and as a consequence the nodes $k_{m+1}$ and $k'_{m+1}$ at which this residual 
supply is depleted will eventually decrease (or move upwards in Figure \ref{phase1}). However, the same argument as in \eqref{orderkm} with
$V_{m+1}$ replaced by $V_{m+1}-\varepsilon$, implies that they remain ordered as $k'_{m+1}\leq k_{m+1}$. Hence, 
we can continue the flow transfers until the full supply $\pi_{m+1}^{m+1}$ is moved from $m+1$ 
to $m$. This completes the first phase, along which $\Delta^\varepsilon_{m+1,n}$ decreases. 

At this point, the supply at $m$ has been increased to $V_m=\pi_m^{m+1}+\pi_{m+1}^{m+1}$, 
and is being shipped to $j=m,\ldots,k_m$ in $z_\varepsilon$ and to $j=m,\ldots,k'_m$ in $z'_\varepsilon$. 
We further observe that 
$$\pi_m^{m+1}+\pi_{m+1}^{m+1}=1-\sum_{i=0}^{m-1}\pi_i^{m+1}\geq 1-\sum_{i=0}^{m-1}\pi_i^{m} = \pi_m^m.$$
We then start the second phase in which we reduce the supply $\pi_m^{m+1}+\pi_{m+1}^{m+1}$
down to  $\pi_m^m$ by transfering some flow $\varepsilon$ from node $m$ to $m-1$. As we proceed
with these flow transfers once again the nodes $k_m$ and $k'_m$ will decrease, but remain
in the order $k'_m\leq k_m$ as it follows from \eqref{orderkm} with $V_m$ replaced by $V_m-\varepsilon$.
As in the first phase, these flow transfers further reduce $\Delta^\varepsilon_{m+1,n}$. The second phase
stops as soon as the supply  $\pi_m^{m+1}+\pi_{m+1}^{m+1}$ has been reduced to $\pi_m^m$. 

At this point, the supply of $m-1$ is $\pi_{m-1}^{m+1}+\pi_m^{m+1}+\pi_{m+1}^{m+1}-\pi_{m}^m$, 
which is larger than $\pi_{m-1}^m$ as it follows from 
$$\pi_{m-1}^{m+1}+\pi_m^{m+1}+\pi_{m+1}^{m+1}=1-\sum_{i=0}^{m-2}\pi_j^{m+1}\geq 1-\sum_{i=0}^{m-2}\pi_j^{m}=\pi^m_{m-1}+\pi_{m}^{m}.$$
We note that the sum of the modified supplies of $m-1$ and $m$ continues to be equal to $V_{m-1}$ and is shipped to 
$m-1,\ldots,k_{m-1}$ in $z_\varepsilon$ and to $m-1,\ldots,k'_{m-1}$ in $z'_\varepsilon$. 
We may then proceed 
as before by moving a flow $\varepsilon$ from $m-1$ to $m-2$, until the supply of $m-1$ is decreased to $\pi_{m-1}^m$.

We continue in this manner with $i=m-2,m-3,\ldots,0$ until the transport plans $z_\varepsilon$ and $z'_\varepsilon$ are 
shipping the correct supplies $\pi_i^m$ from each source node $i=0,\ldots,m$. This is possible since at the end of each
phase the supply of node $i$ is equal to $\sum_{j=i}^{m+1}\pi_j^{m+1}-\sum_{j=i+1}^m\pi_j^m$ 
which is larger than $\pi_i^m$ as it follows as before from the inequality
$$\sum_{j=i}^{m+1}\pi_j^{m+1}=1-\sum_{j=0}^{i-1}\pi_j^{m+1}\geq 1-\sum_{j=0}^{i-1}\pi_j^{m}=\sum_{j=i}^{m}\pi_j^{m}.$$

The final transport plans $z_\varepsilon$ and $z'_\varepsilon$ are therefore feasible for $d_{m,n+1}$ and $d_{m,n}$ 
respectively. Since they are also simple and nested, it follows from the induction hypothesis $(\tilde Q_{n})$
and Proposition \ref{Le1} that these transformed flows are optimal, so 
the final value of $\Delta^\varepsilon_{m+1,n}$ coincides with  $\Delta_{m,n}$ and therefore
$\Delta_{m,n}=\Delta^\varepsilon_{m+1,n}\leq\Delta_{m+1,n}$, completing the induction step.
\end{proof}
 
\vspace{1ex}\noindent {\sc Remark 5.} 
The example in Remark 4 also shows that the convex quadrangle inequality does not hold in general,
not even under stochastic dominance. In that example $(Q)$ fails for several combinations $i<j<k<l$,
for instance with $d_{1,5}+d_{2,4}=\frac{139}{100}>\frac{135}{100}=d_{1,4}+d_{2,5}$.
 
\section{Some final examples and comments}\label{Sexamples}
This section presents some interesting examples that expose the (sometimes) surprisingly 
convoluted structure of the bounds $d_{m,n}$. 

\vspace{1ex}
\begin{example}\label{Example1}
The case where each $\pi^n=\delta^n$ is a Dirac mass at $n \in \NN $, induces the trivial distance $d_{m,n}= \mathbbm{1}_{m\neq n}$ in $\overline \NN$.
\end{example}

\vspace{1ex}
\begin{example}\label{Example2}
 The next simplest case is probably when $\pi^n\!=\! (1\!-\!\alpha) \delta^{n-1}\! +\! \alpha \delta^n$
with $\alpha\!\geq\!\frac{1}{2}$.  
Since $\pi_{n-1}^n=1-\alpha \leq \alpha = \pi_n^n$, condition   {\sc (h)} is satisfied and we may run the 
Inside-Out algorithm from \S \ref{squadrangle} to obtain $d_{m,n}(\alpha)$ as 
polynomials in $\alpha$. Figure \ref{d7n} shows these polynomials for $m=7$.
An interesting feature is that in the range $\lfloor m/2\rfloor\leq n<2(m+1)$ with $m\neq n$ the polynomial $d_{m,n}(\alpha)$ 
has degree $\min\{m,n\}+1$ with some peculiar integer coefficients, while oustide this range $d_{m,n}\equiv 1$. 
Naturally, for $\alpha=1$ we recover $d_{m,n}= \mathbbm{1}_{m\neq n}$ as in the  simple Dirac case of Example \ref{Example1}. 
However, it is amazing how a minor tweak in the $\pi^n$'s induces a drastic change in the 
structure of the recursive bounds $d_{m,n}$.

\vspace{-3ex}
\begin{figure}[ht]
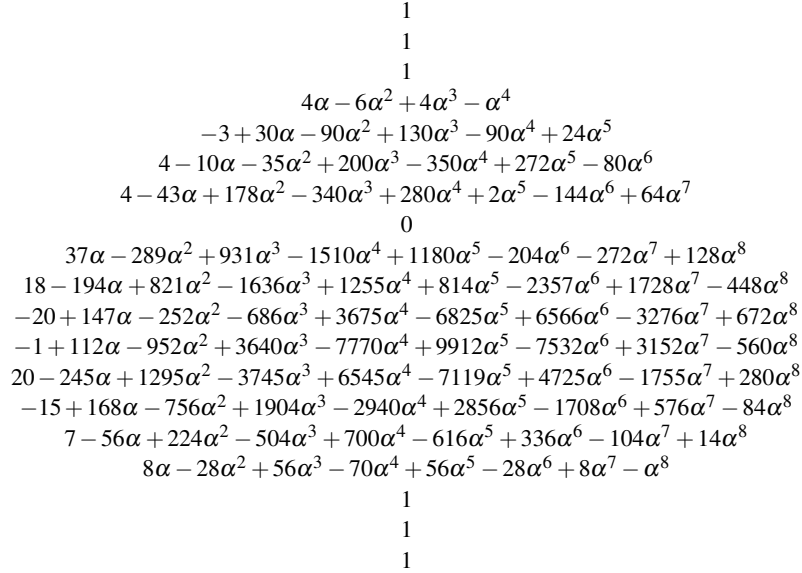

{\scriptsize
$$
\begin{array}{c}
 1 \\
 1 \\
 1 \\
 4 \alpha-6 \alpha ^2+4 \alpha ^3-\alpha ^4  \\
 -3+30 \alpha -90 \alpha ^2+130 \alpha ^3-90 \alpha ^4+ 24 \alpha ^5\\

 4-10 \alpha-35 \alpha ^2+200 \alpha ^3-350 \alpha ^4+272 \alpha ^5 -80 \alpha ^6\\
 4-43 \alpha+178 \alpha ^2-340 \alpha ^3+280 \alpha ^4+2 \alpha ^5-144 \alpha ^6 + 64 \alpha ^7 \\

0\\
37\alpha - 289 \alpha^2 + 931 \alpha^3 - 1510 \alpha^4 +  1180 \alpha^5 - 204 \alpha^6 - 272 \alpha^7 + 128 \alpha^8\\
18 - 194 \alpha + 821 \alpha^2 - 1636 \alpha^3 +  1255 \alpha^4 + 814 \alpha^5 - 2357 \alpha^6 + 1728 \alpha^7 - 448 \alpha^8\\
-20 + 147 \alpha - 252 \alpha^2 - 686 \alpha^3 +  3675 \alpha^4 - 6825 \alpha^5 + 6566 \alpha^6 -  3276 \alpha^7 + 672 \alpha^8\\
-1 + 112 \alpha - 952 \alpha^2 + 3640 \alpha^3 -  7770 \alpha^4 + 9912 \alpha^5 - 7532 \alpha^6 + 3152 \alpha^7 - 560 \alpha^8\\
20 - 245 \alpha + 1295 \alpha^2 - 3745 \alpha^3 +  6545 \alpha^4 - 7119 \alpha^5 + 4725 \alpha^6 - 1755 \alpha^7 + 280 \alpha^8\\
-15 + 168 \alpha - 756 \alpha^2 + 1904 \alpha^3 -  2940 \alpha^4 + 2856 \alpha^5 - 1708 \alpha^6 + 576 \alpha^7 - 84 \alpha^8\\
7 - 56 \alpha + 224 \alpha^2 - 504 \alpha^3 + 700 \alpha^4 - 616 \alpha^5 + 336 \alpha^6 - 104 \alpha^7 + 14 \alpha^8\\
8\alpha - 28 \alpha^2 + 56 \alpha^3 - 70 \alpha^4 +  56 \alpha^5 - 28 \alpha^6 + 8 \alpha^7 - \alpha^8\\
1\\
1\\
1\\
\end{array}
$$
\vspace{-2ex}
}
\caption{\label{d7n}$\{d_{7,n}(\alpha)\}_{0\leq n\leq 18}$ for $\pi^n=(1-\alpha)\delta^{n-1}+\alpha\delta^n$, $\alpha\geq\frac{1}{2}$.}
\end{figure}
\end{example}

\vspace{1ex}
\begin{example}\label{Example3}
 Let us consider next the classical Krasnosel'skii iteration with constant $\alpha_n\equiv\alpha$ and  $\pi^n=(1-\alpha)\pi^{n-1}+\alpha \delta^n$. 
Since {\sc (h)} holds automatically, Theorem \ref{teoint} and \S \ref{s51} imply the tight estimate 
 $\|x^n-Tx^n\|\leq R_n= d_{n,n+1}/\alpha$, which can be computed  by the Inside-Out algorithm.
 For $\alpha\geq\frac{1}{2}$ the algorithm takes a particularly simple form, and yields a polynomial expression
 for $R_n$.  Figure \ref{RnK} shows some of these polynomials as a function of $\alpha\geq \frac{1}{2}$.
\begin{figure}[ht]
{\scriptsize
$$
\begin{array}{rcl}
R_1&=&1 - \alpha + \alpha^2\\[0.55ex]
R_2&=&1 - 2 \alpha + 4 \alpha^2 - 4 \alpha^3 + 2 \alpha^4\\[1ex]
R_3&=&1 - 3 \alpha + 9 \alpha^2 - 18 \alpha^3 + 25 \alpha^4 - 21 \alpha^5 + 9 \alpha^6 - \alpha^7\\[1ex]
R_4&=&1 - 4 \alpha + 16 \alpha^2 - 48 \alpha^3 + 112 \alpha^4 - 192 \alpha^5 + 230 \alpha^6 - 180 \alpha^7\\&& + 84 \alpha^8 - 20 \alpha^9 + 2 \alpha^{10}\\[1ex]
R_5&=&1 - 5 \alpha + 25 \alpha^2 - 100 \alpha^3 + 331 \alpha^4 - 876 \alpha^5 + 1795 \alpha^6 - 2762 \alpha^7\\
&&{} + 3106 \alpha^8 - 2482 \alpha^9 + 1366 \alpha^{10} - 500 \alpha^{11} + 117 \alpha^{12} - 16 \alpha^{13} + \alpha^{14} \\[1ex]
R_6&=&1 - 6 \alpha + 36 \alpha^2 - 180 \alpha^3 + 775 \alpha^4 - 2806 \alpha^5 + 8324 \alpha^6 - 19778 \alpha^7 \\
&&{} + 37023 \alpha^8 - 53948 \alpha^9 + 60623 \alpha^{10}   - 52122 \alpha^{11} + 34044 \alpha^{12} - 16770 \alpha^{13} \\
&&{} + 6163 \alpha^{14} - 1652 \alpha^{15} + 308 \alpha^{16}  - 36 \alpha^{17} + 2 \alpha^{18} \\[1ex]
R_7&=&1 - 7 \alpha + 49 \alpha^2 - 294 \alpha^3 + 1562 \alpha^4 - 7222 \alpha^5 + 28408 \alpha^6 - 93187 \alpha^7\\
  &&{} + 251365 \alpha^8  - 552678 \alpha^9 + 985643 \alpha^{10} - 1422448 \alpha^{11} +   1660135 \alpha^{12} \\
  &&{} - 1567511 \alpha^{13} + 1198337 \alpha^{14} - 741914 \alpha^{15} + 371352 \alpha^{16} - 149443 \alpha^{17} \\
  &&{} + 47802 \alpha^{18} - 11909 \alpha^{19} + 2233 \alpha^{20} - 297 \alpha^{21} + 25 \alpha^{22} - \alpha^{23} 
\end{array}
$$
}
\caption{\label{RnK} $R_n(\alpha)$ for Krasnosel'skii iteration with $\alpha_n\equiv\alpha\geq\frac{1}{2}$.}
\end{figure}
Again, we know very little about these polynomials $R_n(\alpha)$, except that the first two leading terms 
can be proved to be $1-n \alpha$, while the quadratic and cubic terms seem to be $n^2 \alpha^2$ and $(n^2-n^3)\alpha^3$. 
On the other hand, by considering a family of sub-optimal transport plans in $ (\Pcal_{m,n})$, \cite{csv} established
 a non-trivial bound expressed in terms of Gauss' hypergeometric function  
\[
R_n(\alpha) \leq \mbox{${}_2F_{1}\left (-n, \frac{1}{2};2;4\alpha(1\!-\!\alpha) \right )$} \quad \text{for all} \,\, n\geq 1.
\]
Interestingly, the following observation hints at a possible underlying combinatorial structure of the polynomials $R_n$. 
The sequence of degrees of $R_n$ 
$$
\mbox{\small $\{2,4,7,10,14,18,23,28,34,40,47,54,62,70,79,88,98,108,119,130,\ldots\}$}
$$   
coincides with $\lfloor(n^2 + 6 n + 1)/4\rfloor$ at least up to $n=40$. This matches the sequence A014616 in the 
{\em Online Encyclopedia of Integer Sequences} ({\ttfamily oeis.org}), which is the solution  of the 
{\em postage stamp problem} with $n$ stamps and $k=2$ denominations.
The postage stamp problem is as follows: an envelope can carry at most $n$ stamps with values
chosen from $k$ possible denominations $\{ a_1, \ldots, a_k\}$. The problem is 
to find denominations that allow to make up all consecutive postage values 
$p=1,2,\ldots,N$ with $N$ as large as possible. For instance, when $n=3$ and $k=2$ the solution is $N=7$ and is 
attained with $a_1=1, a_2=3$ as follows
{\small
\[
1,1+1,1+1+1,1+3,1+1+3,3+3, 1+3+3.
\]
}
\end{example}

\vspace{1ex}
\begin{example}\label{Example4}
 Another interesting common case is when $\pi^n$ is chosen as a uniform distribution over $\{0, \ldots, n\}$ with $\pi^n_i= \frac{1}{n+1}$ for all $0 \leq i \leq n$. 
This is again a classical  Krasnosel'skii iteration since $\pi_i^n =(1-\alpha_n)\pi^{n-1} + \alpha_n\delta^{n}$ with
$\alpha_n= \frac{1}{n+1}$. Explicitly, in this case the iteration {\sc (km)} reads 
 $$\mbox{$x^{n}=\frac{n}{n+1}x^{n-1} + \frac{1}{n+1}Tx^{n-1}, \quad n \geq 1.$}$$ 
The matrix below presents the exact values of $d_{m,n}$ for $0\leq m,n\leq 6$. 
$$
\arraycolsep=1pt\def\arraystretch{1.4}
d_{m,n}\text{ = }\left(
\begin{array}{ccccccc}
 0 & \frac{1}{2} & \frac{2}{3} & \frac{3}{4} & \frac{4}{5} & \frac{5}{6} & \frac{6}{7} \\
 \frac{1}{2} & 0 & \frac{1}{4} & \frac{3}{8} & \frac{7}{15} & \frac{19}{36} & \frac{97}{168} \\
 \frac{2}{3} & \frac{1}{4} & 0 & \frac{23}{144} & \frac{47}{180} & \frac{1}{3} & \frac{47}{120} \\
 \frac{3}{4} & \frac{3}{8} & \frac{23}{144} & 0 & \frac{329}{2880} & \frac{1681}{8640} & \frac{1733}{6720} \\
 \frac{4}{5} & \frac{7}{15} & \frac{47}{180} & \frac{329}{2880} & 0 & \frac{7609}{86400} & \frac{7793}{50400} \\
 \frac{5}{6} & \frac{19}{36} & \frac{1}{3} & \frac{1681}{8640} & \frac{7609}{86400} & 0 & \frac{257219}{3628800} \\
 \frac{6}{7} & \frac{97}{168} & \frac{47}{120} & \frac{1733}{6720} & \frac{7793}{50400} & \frac{257219}{3628800} & 0 \\
\end{array}
\right)
$$

\vspace{1ex}
\noindent Regarding the bounds $R_n=d_{n,n+1}/\alpha_{n+1}$, from \eqref{bound} we readily obtain
$$R_n=\frac{d_{n,n+1}}{\alpha_{n+1}}\leq\frac{\mathop{\rm diam} (C)}{\sqrt{\pi\sum_{k=1}^n\alpha_k(1\!-\!\alpha_k)}}  \sim  O \left ( \frac{1}{\sqrt{\ln n}} \right ).$$
Now, as mentioned earlier, the tightness of the inequality \eqref{bound} proved in \cite{brc} 
occurs with $\alpha_n\equiv\alpha \approx 1$, and since here $\alpha_n =\frac{1}{n+1}\to 0$ one could 
expect a faster rate for $R_n$. However, numerical calculations suggest that this 
is not the case and the bound $O \left (1/\sqrt{\ln n} \right )$ seems to capture accurately the asymptotic behavior of $R_n$.  Figure~\ref{uniforme} displays $\varphi_n =R_n \sqrt{\ln(n+1)}$ for $1\leq n \leq 500$. For the sake of comparison we also include $\phi_n =R_n\ln(n+1)$.

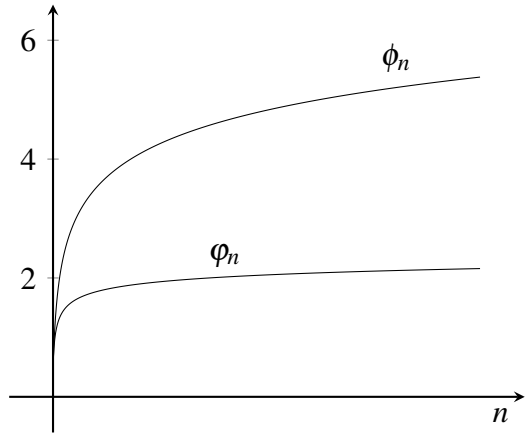
\begin{figure}[h!]
\begin{tikzpicture}[scale=1]
    \begin{axis}[
        table/header=false,
        table/row sep=\\,
        xtick=\empty,
        domain=0:520,
       xmin= -1.0,
      xmax=500,
       ymin= 0,
      ymax=6,
     axis line style = thick,
      axis lines = middle,
      enlargelimits = true,
      ytick={},
    ]
    \node[below right] at (axis cs:500,0) {$n$};
 \node[anchor=south] at (axis cs:400,5.27) {$\phi_n$};
  \node[anchor=south] at (axis cs:200,2) {$\varphi_n$};
    \addplot [mark=none] table[x expr=\coordindex,y index=0]{0.5533742\\0.8156347\\0.9681812\\1.0757479\\1.1525879 
\\1.2154483
\\1.2647739
\\1.3067965
\\1.3424495
\\1.3739294
\\1.4009465
\\1.4258185
\\1.4477209
\\1.4677186
\\1.4860620
\\1.5030220
\\1.5184780
\\1.5331049
\\1.5465344
\\1.5592381
\\1.5711711
\\1.5824392
\\1.5929837
\\1.6031080
\\1.6126719
\\1.6217735
\\1.6304579
\\1.6388036
\\1.6467131
\\1.6543854
\\1.6616886
\\1.6687378
\\1.6755282
\\1.6820625
\\1.6883587
\\1.6944837
\\1.7003731
\\1.7060781
\\1.7116002
\\1.7169830
\\1.7221717
\\1.7272446
\\1.7321501
\\1.7369267
\\1.7415871
\\1.7461208
\\1.7505224
\\1.7548378
\\1.7590331
\\1.7631344
\\1.7671361
\\1.7710530
\\1.7748695
\\1.7786126
\\1.7822679
\\1.7858526
\\1.7893603
\\1.7927966
\\1.7961533
\\1.7994618
\\1.8026957
\\1.8058663
\\1.8089814
\\1.8120395
\\1.8150390
\\1.8179911
\\1.8208851
\\1.8237321
\\1.8265276
\\1.8292825
\\1.8319828
\\1.8346491
\\1.8372670
\\1.8398423
\\1.8423802
\\1.8448784
\\1.8473376
\\1.8497635
\\1.8521479
\\1.8545023
\\1.8568221
\\1.8591085
\\1.8613591
\\1.8635840
\\1.8657769
\\1.8679393
\\1.8700722
\\1.8721797
\\1.8742552
\\1.8763086
\\1.8783339
\\1.8803340
\\1.8823093
\\1.8842596
\\1.8861851
\\1.8880910
\\1.8899712
\\1.8918299
\\1.8936667
\\1.8954842
\\1.8972785
\\1.8990545
\\1.9008087
\\1.9025441
\\1.9042620
\\1.9059602
\\1.9076385
\\1.9093019
\\1.9109453
\\1.9125732
\\1.9141827
\\1.9157778
\\1.9173547
\\1.9189169
\\1.9204631
\\1.9219939
\\1.9235104
\\1.9250115
\\1.9264971
\\1.9279711
\\1.9294298
\\1.9308747
\\1.9323062
\\1.9337247
\\1.9351294
\\1.9365231
\\1.9379028
\\1.9392710
\\1.9406264
\\1.9419709
\\1.9433023
\\1.9446235
\\1.9459331
\\1.9472308
\\1.9485182
\\1.9497951
\\1.9510606
\\1.9523165
\\1.9535609
\\1.9547966
\\1.9560218
\\1.9572370
\\1.9584420
\\1.9596387
\\1.9608256
\\1.9620029
\\1.9631714
\\1.9643311
\\1.9654811
\\1.9666237
\\1.9677567
\\1.9688817
\\1.9699982
\\1.9711067
\\1.9722067
\\1.9732994
\\1.9743838
\\1.9754603
\\1.9765289
\\1.9775907
\\1.9786444
\\1.9796913
\\1.9807305
\\1.9817625
\\1.9827879
\\1.9838060
\\1.9848167
\\1.9858216
\\1.9868191
\\1.9878103
\\1.9887949
\\1.9897731
\\1.9907446
\\1.9917099
\\1.9926690
\\1.9936223
\\1.9945693
\\1.9955103
\\1.9964448
\\1.9973745
\\1.9982977
\\1.9992156
\\2.0001275
\\2.0010340
\\2.0019349
\\2.0028305
\\2.0037206
\\2.0046052
\\2.0054848
\\2.0063593
\\2.0072281
\\2.0080925
\\2.0089515
\\2.0098054
\\2.0106546
\\2.0114990
\\2.0123382
\\2.0131732
\\2.0140029
\\2.0148285
\\2.0156492
\\2.0164653
\\2.0172768
\\2.0180842
\\2.0188870
\\2.0196853
\\2.0204793
\\2.0212692
\\2.0220545
\\2.0228362
\\2.0236133
\\2.0243865
\\2.0251555
\\2.0259205
\\2.0266813
\\2.0274386
\\2.0281918
\\2.0289411
\\2.0296864
\\2.0304282
\\2.0311662
\\2.0319005
\\2.0326308
\\2.0333577
\\2.0340811
\\2.0348008
\\2.0355168
\\2.0362296
\\2.0369386
\\2.0376443
\\2.0383467
\\2.0390457
\\2.0397411
\\2.0404334
\\2.0411224
\\2.0418081
\\2.0424905
\\2.0431698
\\2.0438457
\\2.0445189
\\2.0451887
\\2.0458554
\\2.0465191
\\2.0471797
\\2.0478373
\\2.0484920
\\2.0491437
\\2.0497925
\\2.0504383
\\2.0510814
\\2.0517213
\\2.0523588
\\2.0529933
\\2.0536249
\\2.0542539
\\2.0548802
\\2.0555036
\\2.0561244
\\2.0567425
\\2.0573581
\\2.0579710
\\2.0585812
\\2.0591888
\\2.0597940
\\2.0603965
\\2.0609967
\\2.0615942
\\2.0621893
\\2.0627818
\\2.0633721
\\2.0639598
\\2.0645451
\\2.0651282
\\2.0657087
\\2.0662869
\\2.0668630
\\2.0674366
\\2.0680079
\\2.0685768
\\2.0691437
\\2.0697082
\\2.0702706
\\2.0708307
\\2.0713886
\\2.0719444
\\2.0724980
\\2.0730493
\\2.0735988
\\2.0741460
\\2.0746911
\\2.0752342
\\2.0757752
\\2.0763140
\\2.0768510
\\2.0773859
\\2.0779187
\\2.0784497
\\2.0789786
\\2.0795054
\\2.0800305
\\2.0805536
\\2.0810748
\\2.0815940
\\2.0821113
\\2.0826268
\\2.0831404
\\2.0836521
\\2.0841621
\\2.0846702
\\2.0851765
\\2.0856808
\\2.0861836
\\2.0866845
\\2.0871835
\\2.0876810
\\2.0881766
\\2.0886705
\\2.0891627
\\2.0896531
\\2.0901419
\\2.0906290
\\2.0911144
\\2.0915981
\\2.0920803
\\2.0925608
\\2.0930396
\\2.0935168
\\2.0939924
\\2.0944663
\\2.0949388
\\2.0954097
\\2.0958789
\\2.0963467
\\2.0968128
\\2.0972774
\\2.0977405
\\2.0982021
\\2.0986622
\\2.0991207
\\2.0995778
\\2.1000334
\\2.1004875
\\2.1009402
\\2.1013914
\\2.1018412
\\2.1022895
\\2.1027364
\\2.1031819
\\2.1036259
\\2.1040686
\\2.1045099
\\2.1049498
\\2.1053883
\\2.1058255
\\2.1062613
\\2.1066957
\\2.1071288
\\2.1075606
\\2.1079909
\\2.1084202
\\2.1088480
\\2.1092745
\\2.1096996
\\2.1101236
\\2.1105462
\\2.1109676
\\2.1113877
\\2.1118065
\\2.1122241
\\2.1126405
\\2.1130555
\\2.1134695
\\2.1138821
\\2.1142935
\\2.1147038
\\2.1151128
\\2.1155206
\\2.1159273
\\2.1163327
\\2.1167370
\\2.1171401
\\2.1175421
\\2.1179428
\\2.1183425
\\2.1187410
\\2.1191384
\\2.1195346
\\2.1199297
\\2.1203237
\\2.1207166
\\2.1211083
\\2.1214990
\\2.1218886
\\2.1222771
\\2.1226645
\\2.1230508
\\2.1234361
\\2.1238203
\\2.1242034
\\2.1245855
\\2.1249665
\\2.1253465
\\2.1257255
\\2.1261034
\\2.1264803
\\2.1268562
\\2.1272310
\\2.1276049
\\2.1279777
\\2.1283495
\\2.1287204
\\2.1290902
\\2.1294591
\\2.1298270
\\2.1301939
\\2.1305599
\\2.1309248
\\2.1312889
\\2.1316519
\\2.1320141
\\2.1323752
\\2.1327355
\\2.1330948
\\2.1334531
\\2.1338106
\\2.1341671
\\2.1345227
\\2.1348774
\\2.1352311
\\2.1355840
\\2.1359360
\\2.1362871
\\2.1366373
\\2.1369866
\\2.1373350
\\2.1376825
\\2.1380292
\\2.1383750
\\2.1387199
\\2.1390640
\\2.1394072
\\2.1397496
\\2.1400911
\\2.1404318
\\2.1407716
\\2.1411106
\\2.1414488
\\2.1417861
\\2.1421226
\\2.1424583
\\2.1427932
\\2.1431273
\\2.1434605
\\2.1437930
\\2.1441246
\\2.1444555
\\2.1447855
\\2.1451148
\\2.1454432
\\2.1457710
\\2.1460978
\\2.1464240
\\2.1467494
\\2.1470740
\\2.1473978
\\2.1477209
\\2.1480432
\\2.1483647
\\2.1486855
\\2.1490056
\\2.1493249
\\2.1496435
\\2.1499613
\\2.1502784
\\2.1505948
\\2.1509104
\\2.1512253
\\2.1515395
\\2.1518530
\\2.1521657
\\2.1524778
\\2.1527891
\\2.1530997
\\2.1534097
\\2.1537189
\\2.1540274
\\2.1543353
\\2.1546424
\\2.1549489
\\2.1552546
\\2.1555597
\\2.1558641
\\2.1561679
\\2.1564709
\\2.1567733
\\2.1570750
\\2.1573761
\\2.1576765
\\2.1579762
\\
};

\addplot [mark=none] table[x expr=\coordindex,y index=0]{
  0.4607143
   \\0.8549052
   \\1.1399462
   \\1.3647327
   \\1.5428153
   \\1.6955003
   \\1.8238380
   \\1.9370695
   \\2.0370692
   \\2.1275487
   \\2.2083942
   \\2.2835138
   \\2.3518481
   \\2.4153004
   \\2.4744555
   \\2.5299095
   \\2.5815785
   \\2.6307110
   \\2.6767701
   \\2.7206460
   \\2.7623326
   \\2.8020769
   \\2.8398275
   \\2.8761722
   \\2.9109047
   \\2.9442390
   \\2.9762911
   \\3.0072362
   \\3.0369236
   \\3.0657451
   \\3.0934804
   \\3.1203646
   \\3.1464084
   \\3.1716350
   \\3.1960942
   \\3.2199283
   \\3.2430292
   \\3.2655074
   \\3.2873774
   \\3.3087346
   \\3.3294838
   \\3.3497860
   \\3.3695505
   \\3.3888603
   \\3.4077485
   \\3.4262021
   \\3.4442173
   \\3.4618908
   \\3.4791625
   \\3.4960894
   \\3.5126665
   \\3.5289279
   \\3.5448477
   \\3.5604844
   \\3.5758138
   \\3.5908745
   \\3.6056578
   \\3.6201786
   \\3.6344241
   \\3.6484609
   \\3.6622395
   \\3.6757852
   \\3.6891173
   \\3.7022354
   \\3.7151391
   \\3.7278539
   \\3.7403602
   \\3.7526834
   \\3.7648165
   \\3.7767840
   \\3.7885592
   \\3.8001866
   \\3.8116383
   \\3.8229286
   \\3.8340695
   \\3.8450584
   \\3.8558980
   \\3.8666021
   \\3.8771550
   \\3.8875824
   \\3.8978758
   \\3.9080392
   \\3.9180688
   \\3.9279875
   \\3.9377828
   \\3.9474592
   \\3.9570201
   \\3.9664753
   \\3.9758116
   \\3.9850513
   \\3.9941827
   \\4.0032128
   \\4.0121435
   \\4.0209751
   \\4.0297092
   \\4.0383574
   \\4.0469079
   \\4.0553701
   \\4.0637445
   \\4.0720368
   \\4.0802396
   \\4.0883645
   \\4.0964040
   \\4.1043659
   \\4.1122532
   \\4.1200620
   \\4.1277925
   \\4.1354562
   \\4.1430425
   \\4.1505615
   \\4.1580080
   \\4.1653907
   \\4.1727023
   \\4.1799507
   \\4.1871338
   \\4.1942530
   \\4.2013116
   \\4.2083073
   \\4.2152406
   \\4.2221205
   \\4.2289392
   \\4.2357003
   \\4.2424053
   \\4.2490552
   \\4.2556491
   \\4.2621933
   \\4.2686815
   \\4.2751196
   \\4.2815050
   \\4.2878424
   \\4.2941272
   \\4.3003664
   \\4.3065568
   \\4.3126984
   \\4.3187950
   \\4.3248465
   \\4.3308509
   \\4.3368126
   \\4.3427279
   \\4.3486031
   \\4.3544346
   \\4.3602234
   \\4.3659700
   \\4.3716785
   \\4.3773457
   \\4.3829724
   \\4.3885608
   \\4.3941108
   \\4.3996211
   \\4.4050965
   \\4.4105328
   \\4.4159334
   \\4.4212977
   \\4.4266270
   \\4.4319201
   \\4.4371804
   \\4.4424052
   \\4.4475963
   \\4.4527533
   \\4.4578790
   \\4.4629711
   \\4.4680323
   \\4.4730607
   \\4.4780577
   \\4.4830250
   \\4.4879611
   \\4.4928657
   \\4.4977429
   \\4.5025891
   \\4.5074068
   \\4.5121958
   \\4.5169565
   \\4.5216884
   \\4.5263930
   \\4.5310700
   \\4.5357207
   \\4.5403441
   \\4.5449410
   \\4.5495105
   \\4.5540565
   \\4.5585751
   \\4.5630693
   \\4.5675379
   \\4.5719822
   \\4.5764016
   \\4.5807971
   \\4.5851686
   \\4.5895160
   \\4.5938404
   \\4.5981419
   \\4.6024193
   \\4.6066755
   \\4.6109084
   \\4.6151188
   \\4.6193075
   \\4.6234749
   \\4.6276198
   \\4.6317445
   \\4.6358469
   \\4.6399295
   \\4.6439909
   \\4.6480317
   \\4.6520520
   \\4.6560530
   \\4.6600338
   \\4.6639947
   \\4.6679361
   \\4.6718586
   \\4.6757611
   \\4.6796461
   \\4.6835115
   \\4.6873584
   \\4.6911870
   \\4.6949972
   \\4.6987890
   \\4.7025639
   \\4.7063207
   \\4.7100598
   \\4.7137812
   \\4.7174861
   \\4.7211736
   \\4.7248444
   \\4.7284978
   \\4.7321351
   \\4.7357563
   \\4.7393608
   \\4.7429488
   \\4.7465215
   \\4.7500777
   \\4.7536184
   \\4.7571436
   \\4.7606535
   \\4.7641475
   \\4.7676269
   \\4.7710911
   \\4.7745403
   \\4.7779745
   \\4.7813944
   \\4.7847990
   \\4.7881903
   \\4.7915666
   \\4.7949287
   \\4.7982770
   \\4.8016111
   \\4.8049314
   \\4.8082381
   \\4.8115311
   \\4.8148107
   \\4.8180767
   \\4.8213297
   \\4.8245689
   \\4.8277959
   \\4.8310096
   \\4.8342102
   \\4.8373982
   \\4.8405737
   \\4.8437364
   \\4.8468868
   \\4.8500245
   \\4.8531504
   \\4.8562640
   \\4.8593655
   \\4.8624546
   \\4.8655325
   \\4.8685982
   \\4.8716523
   \\4.8746947
   \\4.8777257
   \\4.8807448
   \\4.8837532
   \\4.8867498
   \\4.8897355
   \\4.8927101
   \\4.8956735
   \\4.8986259
   \\4.9015678
   \\4.9044986
   \\4.9074187
   \\4.9103280
   \\4.9132273
   \\4.9161157
   \\4.9189939
   \\4.9218616
   \\4.9247190
   \\4.9275663
   \\4.9304035
   \\4.9332305
   \\4.9360479
   \\4.9388551
   \\4.9416527
   \\4.9444404
   \\4.9472184
   \\4.9499866
   \\4.9527455
   \\4.9554948
   \\4.9582346
   \\4.9609653
   \\4.9636864
   \\4.9663982
   \\4.9691012
   \\4.9717949
   \\4.9744795
   \\4.9771550
   \\4.9798218
   \\4.9824795
   \\4.9851286
   \\4.9877687
   \\4.9904004
   \\4.9930232
   \\4.9956376
   \\4.9982432
   \\5.0008407
   \\5.0034295
   \\5.0060098
   \\5.0085822
   \\5.0111461
   \\5.0137018
   \\5.0162493
   \\5.0187886
   \\5.0213202
   \\5.0238435
   \\5.0263590
   \\5.0288664
   \\5.0313662
   \\5.0338581
   \\5.0363422
   \\5.0388186
   \\5.0412874
   \\5.0437484
   \\5.0462021
   \\5.0486481
   \\5.0510867
   \\5.0535179
   \\5.0559416
   \\5.0583578
   \\5.0607672
   \\5.0631690
   \\5.0655637
   \\5.0679511
   \\5.0703315
   \\5.0727047
   \\5.0750711
   \\5.0774303
   \\5.0797826
   \\5.0821280
   \\5.0844665
   \\5.0867981
   \\5.0891231
   \\5.0914410
   \\5.0937525
   \\5.0960572
   \\5.0983554
   \\5.1006468
   \\5.1029317
   \\5.1052101
   \\5.1074819
   \\5.1097475
   \\5.1120065
   \\5.1142589
   \\5.1165054
   \\5.1187454
   \\5.1209790
   \\5.1232065
   \\5.1254278
   \\5.1276429
   \\5.1298519
   \\5.1320546
   \\5.1342515
   \\5.1364422
   \\5.1386270
   \\5.1408057
   \\5.1429786
   \\5.1451455
   \\5.1473065
   \\5.1494618
   \\5.1516112
   \\5.1537547
   \\5.1558927
   \\5.1580248
   \\5.1601513
   \\5.1622721
   \\5.1643873
   \\5.1664968
   \\5.1686008
   \\5.1706994
   \\5.1727923
   \\5.1748798
   \\5.1769619
   \\5.1790384
   \\5.1811096
   \\5.1831754
   \\5.1852360
   \\5.1872912
   \\5.1893411
   \\5.1913856
   \\5.1934251
   \\5.1954593
   \\5.1974883
   \\5.1995121
   \\5.2015310
   \\5.2035446
   \\5.2055531
   \\5.2075566
   \\5.2095551
   \\5.2115486
   \\5.2135371
   \\5.2155206
   \\5.2174993
   \\5.2194730
   \\5.2214419
   \\5.2234059
   \\5.2253651
   \\5.2273194
   \\5.2292690
   \\5.2312138
   \\5.2331539
   \\5.2350893
   \\5.2370200
   \\5.2389459
   \\5.2408674
   \\5.2427841
   \\5.2446962
   \\5.2466038
   \\5.2485068
   \\5.2504053
   \\5.2522992
   \\5.2541887
   \\5.2560736
   \\5.2579542
   \\5.2598303
   \\5.2617019
   \\5.2635693
   \\5.2654322
   \\5.2672908
   \\5.2691451
   \\5.2709951
   \\5.2728407
   \\5.2746822
   \\5.2765192
   \\5.2783522
   \\5.2801809
   \\5.2820055
   \\5.2838258
   \\5.2856420
   \\5.2874541
   \\5.2892621
   \\5.2910659
   \\5.2928657
   \\5.2946613
   \\5.2964531
   \\5.2982408
   \\5.3000244
   \\5.3018041
   \\5.3035797
   \\5.3053515
   \\5.3071193
   \\5.3088832
   \\5.3106432
   \\5.3123993
   \\5.3141516
   \\5.3158999
   \\5.3176445
   \\5.3193853
   \\5.3211222
   \\5.3228554
   \\5.3245848
   \\5.3263104
   \\5.3280324
   \\5.3297506
   \\5.3314651
   \\5.3331759
   \\5.3348831
   \\5.3365865
   \\5.3382864
   \\5.3399826
   \\5.3416753
   \\5.3433643
   \\5.3450498
   \\5.3467316
   \\5.3484100
   \\5.3500848
   \\5.3517561
   \\5.3534239
   \\5.3550882
   \\5.3567490
   \\5.3584064
   \\5.3600603
   \\5.3617109
   \\5.3633579
   \\5.3650016
   \\5.3666419
   \\5.3682788
   \\5.3699124
   \\5.3715426
   \\5.3731694
   \\5.3747930
   \\5.3764133
   \\5.3780302
   \\5.3796439\\
   };
  \end{axis}
\end{tikzpicture}
\caption{\label{uniforme}  The quantities $\phi_n$ and $\varphi_n$,  $1\leq n\leq 500$.}
\end{figure}

\end{example}

\vskip 6mm
\noindent{\bf Acknowledgments}

\noindent  
The work of M. Bravo was partially funded by FONDECYT  1191924. 
R. Cominetti gratefully acknowledges the support provided by FONDECYT 1171501.

\end{document}